%% file: main.tex
\documentclass[hidelinks,onefignum,onetabnum]{siamart251216}

\input{ex_shared}

\ifpdf
\hypersetup{
  pdftitle={Equivalent Graphical Conditions for Angle Rigidity and Angle-Based Network Localization},
  pdfauthor={J. HUANG, G. JING, K. ZHAO, and Y. SONG}
}
\fi

\begin{document}

\maketitle

\begin{abstract}
Angle rigidity theory provides a theoretical basis for characterizing the geometric shape of graphs embedded in the plane under angle constraints, and establishing global convergence of angle-based network localization algorithms. The related references, however, are largely restricted to graphs generated by triangulation or bilateration ordering, which are sufficient but not necessary for angle rigidity. This paper addresses this challenge by focusing on establishing the equivalence between angle constraints and bearing constraints. The conditions for cycle compatible frameworks (with possibly non-generic configurations) to be angle fixable are provided. Moreover, it is shown for the first time that Laman graphs equivalently characterize minimally globally angle rigid graphs, implying that a Laman graph is the simplest graph that allows a generic embedding in the plane to be uniquely determined by subtended angles. As an application, a distributed algorithm for angle-based sensor network localization is proposed, which guarantees global convergence under a much milder network topology condition than those in the literature.
\end{abstract}



\section{Introduction}
Determining the geometric shape of a network is a core problem in many multi-agent coordination tasks, including network localization \cite{aspnes2006theory,fang2009sequential, wan2019sensor} and formation control \cite{anderson2008rigid, krick2009stabilisation}. Traditional graph rigidity theory (i.e., distance rigidity theory \cite{asimow1978rigidity}) characterizes when local distance constraints among neighboring vertices uniquely determine a shape, and has seen substantial progress in recent years. To accommodate diverse constraints in practical scenarios \cite{oh2015survey}, distance rigidity theory has been extended to rigidity theories based on bearings \cite{michieletto2021unified,eren2012formation, zhao2015bearing,fang20232}, angles \cite{jing2019angle, chen2020angle, huang2026signed}, and ratio of distance \cite{cao2019ratio}. As a result, the conditions to determine a unique shape have been established for different types of local geometric constraints.

Global distance rigidity guarantees shape uniqueness and has been favored in sensor network localization \cite{aspnes2006theory}. The equivalent graphical condition\footnote{Throughout this paper, the “equivalent graphical condition” for some rigidity property refers to the mildest condition on the graph $\mathcal{G}$ for the framework $(\mathcal{G}, p)$ to have this rigidity property. More specifically, there exists a configuration $p$ such that $(\mathcal{G},p)$ has the rigidity property if and only if $\mathcal{G}$ satisfies this graphical condition.} for global distance rigidity in the plane was established in \cite{jackson2005connected}. However, it usually requires a large number of edges. In \cite{zhao2015bearing,zhao2017laman}, a class of sparser graphs (i.e., Laman graphs) has been shown to be the equivalent graphical condition for infinitesimal bearing rigidity, which guarantees shape uniqueness under bearing constraints. However, the bearing-based methods in \cite{zhao2015bearing,zhao2016localizability} require measurements to be expressed in the same coordinate frame. Although the formulation in \cite{michieletto2021unified} accommodates body-frame bearing measurements, the equivalent graphical condition for this type of rigidity remains unclear.

Since angle measurements are independent of the coordinate frame \cite{peng2025angle, jing2021angle, chen2022triangular, jing2018weak}, and an angle-constrained shape admits higher degrees of freedom than distance- and bearing-based approaches \cite{jing2019angle, chen2020angle, huang2026signed}, angle rigidity theory has attracted increasing attention recently. In our previous work \cite{huang2026signed}, we established an equivalent graphical condition for signed angle rigidity, ensuring shape uniqueness under signed angle constraints (taking values in $[0,2\pi)$). However, signed angle measurements inherently assume that all agents share the same convention for clockwise and counterclockwise directions. Unsigned angle measurements\footnote{In the remainder of this paper, ``unsigned angle" or ``angle" refers to the cosine of the angle between edges; the angle rigidity properties refer to those formulated in the rigidity theory of \cite{jing2019angle}.} (taking values in $[0,\pi]$) are invariant under changes in the handedness of local coordinate frames and have been considered in the angle rigidity theory \cite{jing2019angle}. Unfortunately, the related equivalent graphical condition remains to be established. The algorithms developed in \cite{jing2019angle,jing2021angle} for angle-based network localization and formation control are applicable only to graphs generated by triangulation or bilateration ordering. 

In this paper, we investigate graphical conditions for angle rigidity by relating it to the well-established bearing rigidity theory and exploiting several intuitive yet often overlooked structural properties of rigid graphs. This approach is motivated by related developments in graph rigidity theory. The study in \cite{liu2025rigidity} characterized rigidity properties from a cycle perspective, highlighting the role of the cycle structure in the analysis of rigidity.  In \cite{zhao2017laman}, the equivalent graphical condition for infinitesimal bearing rigidity in the plane was derived using its equivalence with infinitesimal distance rigidity. In \cite{huang2026signed}, an equivalence between signed angle rigidity and bearing rigidity was established. Note that unsigned angle measurements carry less information than signed angle measurements (see the discussion below Lemma \ref{theorem: having same R for positively dependent}). The equivalence between angle rigidity and bearing rigidity can not be established by trivially extending the proof approach used in \cite{huang2026signed}.

The main contributions of this paper are summarized as follows. 
\begin{itemize}
    \item[(\rmnum{1})] The equivalence between angle constraints and bearing constraints (up to uniform orthogonal transformations) is established from a cycle perspective; see Lemmas \ref{lemma: cycle fixable can determine the bearing} and \ref{lemma: angle can determine the bearing}.
    \item[(\rmnum{2})] We establish the equivalence between (infinitesimal, global) angle rigidity and (infinitesimal, global) bearing rigidity for cycle compatible frameworks and graphs, and show that Laman graphs equivalently characterize minimally globally angle rigid graphs; see Theorems \ref{theorem: equivalence in framework2} and \ref{Corollary: infinitesimal angle rigidity inplies global rigdity}.
    \item[(\rmnum{3})] As an application, we design a distributed protocol with global convergence guarantee for angle-based sensor network localization under a milder condition on the underlying graph than those in the literature; see Theorem \ref{theorem:localization theorem}.
\end{itemize}

The paper is organized as follows. Section \ref{section:2} formulates the network localization problem and reviews angle rigidity theory. Section \ref{sec:3} investigates the conditions for uniquely determining a shape by angle constraints. Section \ref{sec:4} demonstrates the application of our results to network localization. Section \ref{sec:5} presents a numerical example for the developed algorithm. Section \ref{sec:6} summarizes the entire paper.

\textit{Notations:} In this paper, ${\mathbb{R}}$ denotes the set of real numbers; ${\mathbb{R}}^{d}$ is the $d$--dimensional Euclidean space; $\|\cdot\|$ is the $L_2$ norm; $\otimes$ is the Kronecker product; $|A|$ denotes the cardinality of the set $A$; $X^{\top}$ means the transpose of the matrix $X$; rank($X$) denotes the rank of the matrix $X$; Define $I_n \in \mathbb{R}^{n\times n}$ as the identity matrix; $A\backslash B$ denotes the set of elements in $A$ but not in $B$; $\mathscr{R}_o(\alpha)$ is the 2-dimensional rotation matrix with $\alpha \in [0,2\pi)$; for $x\in {\mathbb{R}}^{2}$, $x^\perp \triangleq \mathscr{R}_o(\frac{\pi}{2})x$; for $X_{i}\in {\mathbb{R}}^{a\times b}, i = 1,\cdots,q$, diag$\left(X_i\right) \triangleq $ blockdiag$\left\{X_1,\cdots,X_q\right\}\in{\mathbb{R}}^{qa\times qb}$. 

Throughout this paper, we consider the undirected graph ${\mathcal{G}} \triangleq \left(\mathcal{V}, \mathcal{E}\right)$, where $\mathcal{V} \triangleq \left\{1,\cdots,n\right\}$ is the vertex set with $n\geq3$, $\mathcal{E} \subseteq \mathcal{V} \times \mathcal{V}$ is the edge set with $m = |\mathcal{E}|$. Here, $(i,j)\in\mathcal{E}$ and $(j,i)\in\mathcal{E}$ are identified as the same edge. The neighbor set of vertex $j$ in $\mathcal{G}$ is denoted by $\mathcal{N}^j_{\mathcal{G}} \triangleq \left\{i\in\mathcal{V}:(j,i)\in\mathcal{E}\right\}$. If $\mathcal{G}$ is clear from the context, we simply denote $\mathcal{N}^j_{\mathcal{G}}$ by $\mathcal{N}_j$. The complete graph with the same vertex set as $\mathcal{G}$ is denoted by $\mathcal{K}$. A graph $\bar{\mathcal{G}} \triangleq (\bar{\mathcal{V}},\bar{\mathcal{E}})$ is said to be a subgraph of $\mathcal{G}$ if $\bar{\mathcal{V}} \subseteq \mathcal{V}$ and $\bar{\mathcal{E}} \subseteq (\bar{\mathcal{V}} \times \bar{\mathcal{V}}) \bigcap \mathcal{E}$. Furthermore, if $\bar{\mathcal{G}}$ is connected and satisfies $\bar{\mathcal{V}} = \mathcal{V}$, it is said to be a spanning subgraph. Within an undirected graph $\mathcal{G}$, an undirected path between $l_1\in\mathcal{V}$ and $l_{n_P}\in\mathcal{V}$ is a sequence of edges of the form $(l_1,l_2)$, $(l_2,l_3)$, $\dots, (l_{n_P-1},l_{n_P})$, in which $l_i \neq l_j$ for $j\in\{1,\cdots, n_P\}\backslash\{i\}$. We use the subgraph $\mathcal{P} \triangleq (\mathcal{V}_P,\mathcal{E}_P)$ to represent such a path, where $\mathcal{V}_P \triangleq \{l_1,\dots,l_{n_P}\} \subseteq \mathcal{V}$ and $\mathcal{E}_P \triangleq \{(l_1,l_2),\dots,(l_{n_P-1},l_{n_P})\} \subseteq\mathcal{E}$. An undirected cycle is a path whose first and last vertices are the same. For subsequent analysis, we represent it by the subgraph $\mathcal{C}_{n_P} \triangleq(\mathcal{V}^{n_P}_C,\mathcal{E}^{n_P}_C)$ where $\mathcal{V}^{n_P}_C \triangleq \{l_1,\dots,l_{n_P}\} \subseteq \mathcal{V}$ and $|\mathcal{V}^{n_P}_C| \geq 3$ and $\mathcal{E}^{n_P}_C \triangleq \{(l_1,l_2),\dots,(l_{n_P-1},l_{n_P}),(l_{n_P},l_1)\} \subseteq\mathcal{E}$. 

\section{Preliminaries} \label{section:2} In this section, we formulate the angle-based sensor network localization problem and show its connection with angle rigidity.
\subsection{Problem formulation}

Consider a sensor network in the plane with a sensor set $\mathcal{V}=\mathcal{A}\bigcup\mathcal{F}$, where $\mathcal{A}\subseteq\mathcal{V}$ denotes the set of anchors with known locations, and $\mathcal{F}=\mathcal{V}\setminus\mathcal{A}$ denotes the set of followers whose locations are to be estimated. Let $p_j = [x_j,y_j]^\top\in\mathbb{R}^2$ be the real location of sensor $j\in\mathcal{V}$ and $p=[p_1^\top,\ldots,p_n^\top]^\top$ collects the real locations of all sensors in the network, where $p_j\neq p_i$ for all distinct $j,i\in\mathcal{V}$. 

\begin{figure}
		\centering
		\includegraphics[width=1\linewidth]{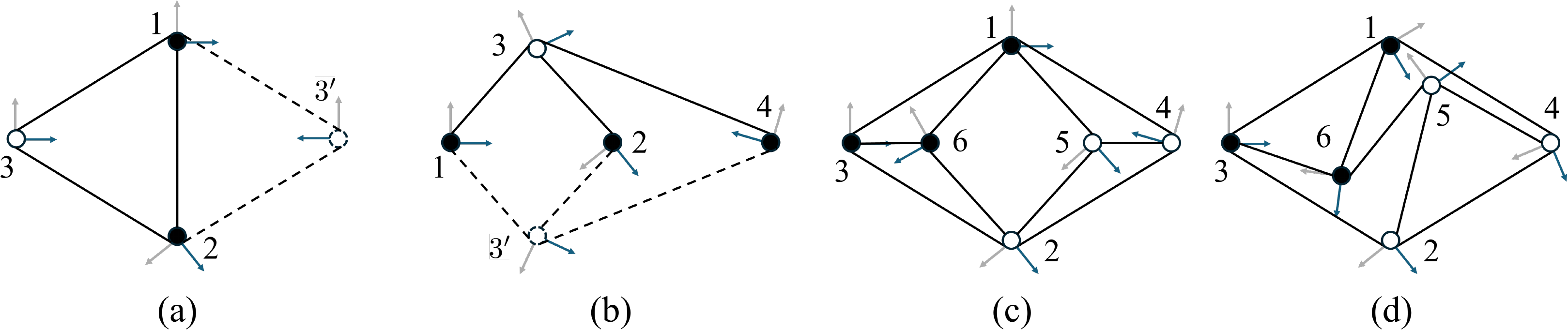} 
		\caption{Sensor networks where black and white vertices are anchors and followers, respectively. The local coordinate frame of each sensor is indicated by the arrows attached to the vertex. The (dash) lines and (dash) vertices in (a) and (b) indicate the possible locations of the network with the same angle measurements and anchors. The edges $(2,3)$ and $(1,4)$ are parallel in (c) and (d).}
		\label{Sensor_network_illustration}
		\centering
	\end{figure}

In this paper, we describe the interaction relationships among sensors by a fixed undirected graph $\mathcal{G}=(\mathcal{V},\mathcal{E})$ and the sensor network by the triplet $(\mathcal{G},p,\mathcal{A})$. In the network, each sensor $j\in\mathcal{V}$ can communicate with its neighboring sensors in $\mathcal{N}^j_{\mathcal{G}}$ to exchange available information. Moreover, sensor $j$ can measure the relative bearings $b^{(j)}_{ji} = (p^{(j)}_{i}-p^{(j)}_{j})/\|p^{(j)}_{i}-p^{(j)}_{j}\|$ pointing from sensor $j$ toward neighbor $i \in \mathcal{N}^j_{\mathcal{G}}$, where $p^{(j)}_{i}$ is the position of sensor $i$ in the local coordinate frame of sensor $j$. Since local frames are inconsistent across sensors, different sensors may have different understandings on clockwise and counterclockwise directions\footnote{Existing localization strategies that rely on a consistent convention on the local frames to obtain measurement information cannot be directly applied in this setting \cite{huang2026signed,zhao2016localizability, oh2026subtended, chen2022triangular}.}. Note that the angle measurements $\{{\rm cos} \theta_{ijk}(p) = (b^{(j)}_{ji})^\top b^{(j)}_{jk} : (i,j,k) \in \mathcal{T}_{\mathcal{G}}\}$ are invariant under orthogonal
changes of the local frames, where $\mathcal{T}_{\mathcal{G}} = \{(i,j,k):(j,i),(j,k)\in\mathcal{E},i < k\}$ indexes all measurable angles. For example, the angle $(b^{(1)}_{12})^\top b^{(1)}_{13}$ measured by sensor 1 in Fig.~\ref{Sensor_network_illustration}(a) has the same value when evaluated in different coordinate frames. From a centralized perspective, the angle-based sensor network localization problem is formulated as the following nonlinear optimization problem:
\begin{eqnarray}\label{equation:optimization problem}
    \mathop{\mathrm{minimize}}\limits_{\hat{p} \in \mathbb{R}^{2n}}&& \sum\limits_{(i,j,k) \in \mathcal{T}_{\mathcal{G}}}(b^\top_{ji}(\hat{p})b_{jk}(\hat{p}) - {\rm cos} \theta_{ijk}(p))^2 \nonumber\\
    \text{subject to}&&\quad \quad \quad \quad \hat{p}_j=p_j, ~~~~j\in\mathcal{A},
\end{eqnarray}
where $\hat{p} = [\hat{p}^\top_1,\dots,\hat{p}^\top_n]^\top$ with $\hat{p}_j\in\mathbb{R}^2$ denoting the estimation of the real location of sensor $j$ for $j\in\mathcal{V}$, $b_{ji}(\hat{p}) = (\hat{p}_i-\hat{p}_j)/\|\hat{p}_i-\hat{p}_j\|$, and ${\rm cos} \theta_{ijk}(p)$ is the angle measurement between the edges $(j,i), (j,k)\in\mathcal{E}$.

Given a set of angle measurements collected from a network, together with the real locations of the anchors, the optimization (\ref{equation:optimization problem}) may have multiple optimal solutions, for each of which the objective function is zero. For example, under the sensor networks shown in Fig. \ref{Sensor_network_illustration}(a) and \ref{Sensor_network_illustration}(b), the global minimizer of (\ref{equation:optimization problem}) is not unique. Therefore, we give the following definition to characterize the solvability of the localization problem.
\begin{definition}
    A sensor network $(\mathcal{G},p,\mathcal{A})$ is said to be angle localizable if the nonlinear least-squares problem (\ref{equation:optimization problem}) has a unique global minimizer.
\end{definition}

Aiming to effectively localize the real locations based on the angle measurements and the known locations, we first propose and address the following problem:
\begin{problem}\label{problem: 1}
    Determine the conditions for $(\mathcal{G},p,\mathcal{A})$ to be angle localizable.
\end{problem}

Once Problem \ref{problem: 1} is solved, the remaining part is to seek the global minimizer of (\ref{equation:optimization problem}). However, solving the optimization (\ref{equation:optimization problem}) directly requires network-wide aggregation of sensor measurements. In this paper, we aim to develop a distributed localization method for each follower $j\in\mathcal{F}$ to estimate its real location $p_j$ using only the local angle measurements $\{{\rm cos}\theta_{ijk}\}_{(i,j,k)\in\mathcal{T}^j_{\mathcal{G}}}$ and information communicated from its neighbors, where $\mathcal{T}^j_{\mathcal{G}} \subseteq \mathcal{T}_{\mathcal{G}}$ is the subset indexing the angles measured at sensor $j$ and its neighbors. We formulate this as the following controller design problem:
\begin{problem} \label{problem: 2}
    Consider an angle localizable sensor network $(\mathcal{G},p,\mathcal{A})$ and an estimator for each sensor: 
    \begin{equation}\label{equation}
        \dot{\hat{p}}_j(t) =  u_j(t), ~~~~j \in \mathcal{V},
    \end{equation}
    design $u_j(t)$ based on $\{{\rm cos} \theta_{ijk}\}_{(i,j,k)\in\mathcal{T}^j_{\mathcal{G}}}$ and information from neighbors such that
\begin{equation}
    \lim\limits_{t\rightarrow \infty}(\hat{p}_j(t)-p_j) = 0, ~~~~j \in \mathcal{V},
\end{equation}
where $\hat{p}_j(t)$ is the estimation of sensor $j$ at time instant $t$ and $\hat{p}_j(t) = p_j$ if $j \in \mathcal{A}$.
\end{problem}
        
In the literature of sensor network localization, the localizability problem is usually studied via graph rigidity theory \cite{chen2022triangular,zhao2016localizability,jing2021angle,huang2026signed}, which answers the question of whether the geometric shape of a framework can be uniquely determined by the local constraints (e.g., angles and bearings) up to rigid transformations (e.g., uniform translations, scalings, or rotations). In order to solve Problem \ref{problem: 1}, in the next subsection, we introduce angle rigidity theory.
    
\subsection{Angle rigidity theory}\label{sec:2} In this subsection, aiming to solve Problem \ref{problem: 1}, we present some definitions and results from angle rigidity theory \cite{jing2019angle}. Detailed tutorials on bearing rigidity theory \cite{zhao2015bearing} are provided in Appendix \ref{appendix: Bearing Rigidity Theory}.

The pair $(\mathcal{G}, p)$ is referred to as a framework in $\mathbb{R}^2$, where $p = \left[p^{\top}_1, \cdots, p^{\top}_n\right]^{\top}$, $p_j = [x_j,y_j]^\top \in \mathbb{R}^2$ is the coordinate mapped to vertex $j\in\mathcal{V}$ with the assumption that $p_i \neq p_j$ for any $i\in\mathcal{V}\backslash\{j\}$. If $p_1,\cdots, p_n$ belong to a hyperplane of $\mathbb{R}^{2}$,  we say that $p$ is \textit{degenerate}, and the framework $(\mathcal{G}, p)$ is said to be a \textit{degenerate framework}. Given a (spanning) subgraph $\bar{\mathcal{G}}$ of $\mathcal{G}$, $(\bar{\mathcal{G}}, \bar{p})$ is said to be the \textit{(spanning) subframework} of $(\mathcal{G},p)$, where $\bar{p}$ is the stack of the coordinates of the vertices in $(\mathcal{G},p)$.

For a framework $(\mathcal{G},p)$ in $\mathbb{R}^2$, the angle constraints can be stacked to form the \textit{angle rigidity function}
\begin{eqnarray} \label{equation: angle rigidity function}
        A_{\mathcal{G}}(p) = [\cdots,b^\top_{ji}b_{jk},\cdots]^\top, ~~~~(i,j,k) \in \mathcal{T}_{\mathcal{G}}, 
\end{eqnarray}
where $b_{ji} \triangleq b_{ji}(p) = (p_i-p_j)/\|p_i-p_j\|\in\mathbb{R}^2$ describes the bearing pointing from $p_j$ to $p_i$ in the global frame and $b^\top_{ji}b_{jk}$ is the cosine of the angle $\theta_{ijk} \triangleq \theta_{ijk}(p)\in[0,\pi]$ between $b_{ji}$ and $b_{jk}$. Denote the signed angle $\alpha_{ji} \triangleq \alpha_{ji}(p) \in [0,2\pi)$ from the $x$-axis (of the global frame) to $b_{ji}$ in the counter-clockwise direction as the angle of $b_{ji}$ \cite{huang2026signed}. 

    By the definition of the angle rigidity function (\ref{equation: angle rigidity function}), the configuration set of frameworks with the same angles as $(\mathcal{G}, p)$ can be represented as
    \begin{eqnarray} \label{equation: configurations of angle constraints}
        A^{-1}_{\mathcal{G}}(A_{\mathcal{G}}(p)) = \{q\in\mathbb{R}^{2n}~:~b^\top_{ji}(q)b_{jk}(q) = b^\top_{ji}(p)b_{jk}(p), (i,j,k)\in\mathcal{T}_{\mathcal{G}}\}.
    \end{eqnarray}
    It can be verified that\footnote{$B^{-1}_{\mathcal{G}}(B_{\mathcal{G}}(p))$ is the configuration set of frameworks with the same bearings as $(\mathcal{G}, p)$. See Appendix \ref{appendix: Bearing Rigidity Theory} for details.} $B^{-1}_{\mathcal{G}}(B_{\mathcal{G}}(p)) \subseteq A^{-1}_{\mathcal{G}}(A_{\mathcal{G}}(p))$. Moreover, since uniform translations, rotations, scalings, and reflections are angle-preserving motions for any framework, we define $\mathscr{E}_A(p) \subseteq A^{-1}_{\mathcal{G}}(A_{\mathcal{G}}(p))$ as
    \begin{eqnarray} 
        \mathscr{E}_A(p) = \{q\in\mathbb{R}^{2n}:q = c(I_n\otimes \mathscr{R})p+\textbf{1}_n \otimes \xi,c\in \mathbb{R}^+, \mathscr{R}\in{\rm O}(2), \xi \in \mathbb{R}^{2} \},
    \end{eqnarray}
    where $\textbf{1}_n \triangleq [1,\cdots,1]^{\top}$. If $A^{-1}_{\mathcal{G}}(A_{\mathcal{G}}(p)) \subseteq \mathscr{E}_A(p)$, then the framework $(\mathcal{G},q)$ with $q\in A^{-1}_{\mathcal{G}}(A_{\mathcal{G}}(p))$ has the same shape as $(\mathcal{G},p)$ up to a similarity transformation.

    \begin{definition}[Angle fixability \cite{jing2019angle}]
        In $\mathbb{R}^2$, a framework $(\mathcal{G}, p)$ is said to be \textit{angle fixable} if $A^{-1}_{\mathcal{G}}(A_{\mathcal{G}}(p)) = \mathscr{E}_A(p)$.
    \end{definition}
    
    Given an angle fixable framework $(\mathcal{G}, p)$ and $q \in A^{-1}_{\mathcal{G}}(A_{\mathcal{G}}(p))$, one has $q=p$ if at least three vertices with $q_i=p_i$ are not collinear \cite{jing2021angle}. For example, the frameworks in Figs. \ref{Sensor_network_illustration}(c) and \ref{Sensor_network_illustration}(d) satisfy $A^{-1}_{\mathcal{G}}(A_{\mathcal{G}}(p)) = \mathscr{E}_A(p)$; together with the known locations of anchors, the configuration $p$ can be determined uniquely. Therefore, angle localizability of a sensor network $(\mathcal{G},p, \mathcal{A})$ is mainly determined by the pair $(\mathcal{G},p)$, and Problem \ref{problem: 1} can be converted into:
    \begin{center}
        \textit{Determine the conditions for $(\mathcal{G},p)$ to be angle fixable.}
    \end{center}
    To answer this question, we provide the definitions associated with angle rigidity.

    \begin{definition}[Angle rigidity and global angle rigidity]
        In $\mathbb{R}^2$, a framework $(\mathcal{G}, p)$ is said to be \textit{angle rigid} if there is a neighborhood $U_p$ of $p$ such that $U_p \bigcap A^{-1}_{\mathcal{G}}(A_{\mathcal{G}} \\ (p)) = U_p \bigcap A^{-1}_{\mathcal{K}}(A_{\mathcal{K}}(p))$, and to be \textit{globally angle rigid} if $A^{-1}_{\mathcal{G}}(A_{\mathcal{G}}(p)) = A^{-1}_{\mathcal{K}}(A_{\mathcal{K}}(p))$.
    \end{definition}

	\begin{figure}
		\centering
		\includegraphics[width=0.93\linewidth]{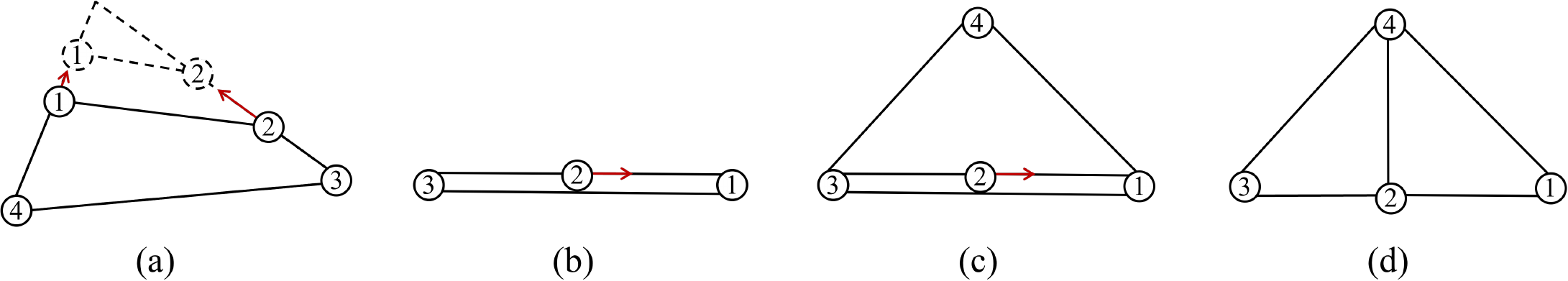} 
		\caption{Illustration of frameworks in which red arrows indicate the non-uniform motions preserving the invariance of angles while deforming the shape. (a) and (c) are not angle rigid. (b) and (d) are both globally angle rigid, but the shape of (b) can not be determined by angles uniquely.}
		\label{example_of_rigidity}
		\centering
	\end{figure}
    
    Fig. \ref{example_of_rigidity} shows several frameworks with different angle rigidity properties. In particular, the framework shown in Fig. \ref{example_of_rigidity}(b) demonstrates that global angle rigidity is not necessary to guarantee $A^{-1}_{\mathcal{G}}(A_{\mathcal{G}}(p)) = \mathscr{E}_A(p)$ due to the existence of continuous non-uniform motions preserving angles but deforming the shape. 
    
    In fact, infinitesimal angle rigidity can circumvent such undesired motions. An instantaneously continuous motion $\dot{p} \triangleq v = [v^\top_1,\dots,v^\top_{n}]^\top$ is said to be an \textit{infinitesimal angle motion} if it satisfies $\dot{A}_{\mathcal{G}}(p) = R_A(p)v = 0 $, where $R_A(p) \triangleq \frac{\partial A_{\mathcal{G}}(p)}{\partial p}$ is the \textit{angle rigidity matrix}. An infinitesimal angle motion is said to be \textit{trivial} if it consists of uniform translations, rotations, and scalings of the entire framework. 

    \begin{definition} [Infinitesimal angle rigidity]
        In $\mathbb{R}^2$, a framework $(\mathcal{G}, p)$ is said to be \textit{infinitesimally angle rigid} if all of its infinitesimal angle motions are trivial.
    \end{definition}

    \begin{remark}
        In what follows, we use shorthands AR, GAR, and IAR to denote angle rigidity (angle rigid), global angle rigidity (globally angle rigid), and infinitesimal angle rigidity (infinitesimally angle rigid), respectively. Similarly, BR, GBR, and IBR are shorthands for corresponding concepts in bearing rigidity theory.
    \end{remark}

The continuous motions denoted by red arrows in Figs. \ref{example_of_rigidity}(a)--\ref{example_of_rigidity}(c) are infinitesimal angle motions but not trivial. Hence, the frameworks are not IAR. In contrast, the framework in Fig. \ref{example_of_rigidity}(d) is IAR. Since it is also GAR, it follows from \cite[Theorem 1]{jing2019angle} that $A^{-1}_{\mathcal{G}}(A_{\mathcal{G}}(p)) = A^{-1}_{\mathcal{K}}(A_{\mathcal{K}}(p)) = \mathscr{E}_A(p)$.
    
In \cite{jing2021angle}, it has been formally proved that a framework in $\mathbb{R}^2$ is angle fixable if and only if it is both IAR and GAR. However, the authors only found a special class of angle fixable frameworks by considering non-degenerate bilateration ordering. In \cite{jing2019angle}, the authors established that a framework embedded by a graph is IAR for almost all configurations as long as it is IAR for one configuration. Although a similar discussion on AR (GAR) has been reported in \cite{jing2019angle}, a rigorous derivation is yet to be provided.  If AR (GAR) is a generic property of the graph, then angle fixability, i.e., the equation $A^{-1}_{\mathcal{G}}(A_{\mathcal{G}}(p)) = \mathscr{E}_A(p)$, is critically determined by the graph.

	\begin{figure}
		\centering
		\includegraphics[width=0.85\linewidth]{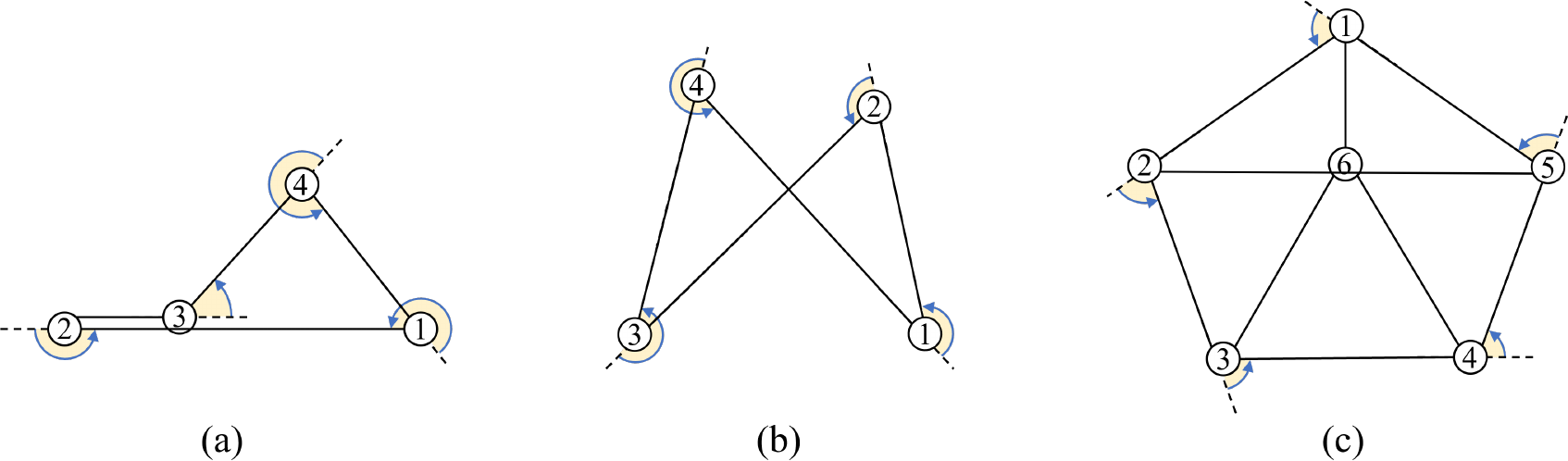} 
		\caption{Illustration of frameworks where the counter-clockwise arc arrows indicate the signed angles. (a) A concave cycle framework $(\mathcal{C}_4,p)$ in which $b_{12}$ is parallel to $b_{23}$ and vertex 3 lies on the edge $(2,1)$. (b) A non-planar cycle framework in which the open edges $(4,1)$ and $(3,2)$ are intersecting. (c) A non-planar framework where the neighbors of vertex 6 are separated by the open edge $(5,2)$ into two parts, the subframework $(\mathcal{C}_5,\bar{p})$ with $\mathcal{V}^{5}_C = \{1,2,3,4,5\}$ is convex.}
		\label{illustration_of_Lemma_proof}
		\centering
	\end{figure}

\section{Angle fixable frameworks and graphs} \label{sec:3} In this section, we focus on Problem \ref{problem: 1} and seek the conditions for guaranteeing angle fixability. By investigating the relationship between angle constraints and bearing constraints, we first show the angle fixability of frameworks with the help of cycle compatibility, then show that the Laman graph characterizes minimally angle rigid graphs.

\subsection{Angle fixability for cycle compatible frameworks}\label{Angle fixability for cycle compatible frameworks} We first show the relationship between angles and bearings in the context of planar cycle frameworks.

In this paper, a framework embedded by a cycle $\mathcal{C}_n =(\mathcal{V}^{n}_C,\mathcal{E}^{n}_C)$ is said to be a \textit{cycle framework}, which has been studied in \cite{arrigoni2018bearing,ko2020bearing, trinh2018formations}. For simplicity, we denote $\mathcal{V}^n_C = \{1, \dots, n\}$ and $\mathcal{E}^n_C = \{(1,2), \dots, (n,1)\}$. Referring to the definition of planar graphs in graph theory \cite[pp. 243]{murty2008graph}, we provide the definition of planar frameworks.

    \begin{definition}[Planar framework] \label{definition plane framework}
        In $\mathbb{R}^2$, $(\mathcal{G},p)$ is said to be planar if no path in the framework is split by any open edge into two components lying in two open half-planes, where an open edge is a line segment with its endpoints removed.
    \end{definition}

    In graph theory, a graph is planar if it admits an embedding in the plane such that each edge is represented by a simple curve and any pair of curves intersect only at common endpoints. In contrast, Definition \ref{definition plane framework} concerns the straight line realization of graphs in the plane. In particular, it allows distinct edges to overlap, while excluding transverse intersections between open edges (see Figs. \ref{illustration_of_Lemma_proof}(a) and \ref{illustration_of_Lemma_proof}(b)). Moreover, it excludes the case where a vertex lying on an open edge  has its neighbors located in the two opposite open half-planes (see Figs. \ref{illustration_of_Lemma_proof}(a) and \ref{illustration_of_Lemma_proof}(c)).

    \begin{definition} \label{definition convex planar cycle framework}
        In $\mathbb{R}^2$, a planar cycle framework $(\mathcal{C}_n,p)$ is said to be convex if $b^\top_{i+1,i+2}b^\bot_{i,i+1}$ has the same sign for all $i\in \mathcal{V}^n_C$, and concave otherwise, where $b_{i,i+1}$ is the bearing from $i$ to $i+1$ (modulo $n$). A convex (or concave) cycle framework $(\mathcal{C}_n,p)$ is said to be strictly convex (or strictly concave) if $b^\top_{i+1,i+2}b^\bot_{i,i+1} \neq 0$ for all $i\in \mathcal{V}^n_C$.
    \end{definition}
    
    Recalling the definition of signed angles in \cite{huang2026signed}, convexity of a planar cycle framework indicates that the signed angles (rotation angles) from $b_{i,i+1}$ to $b_{i+1,i+2}$ in the counter-clockwise direction for all $i\in \mathcal{V}^n_C$ belong to the same range $[0,\pi]$ or $[\pi, 2\pi)$ (see Figs. \ref{illustration_of_Lemma_proof}(a) and \ref{illustration_of_Lemma_proof}(c)), and the sum of these signed angles is $2\pi$ (w.r.t. the range $[0,\pi]$) or $2(n-1)\pi$ (w.r.t. the range $[\pi, 2\pi)$). The following results establish a connection between angle constraints and bearing constraints for cycle frameworks.

 	\begin{figure}[t]
		\centering
		\includegraphics[width=0.75\linewidth]{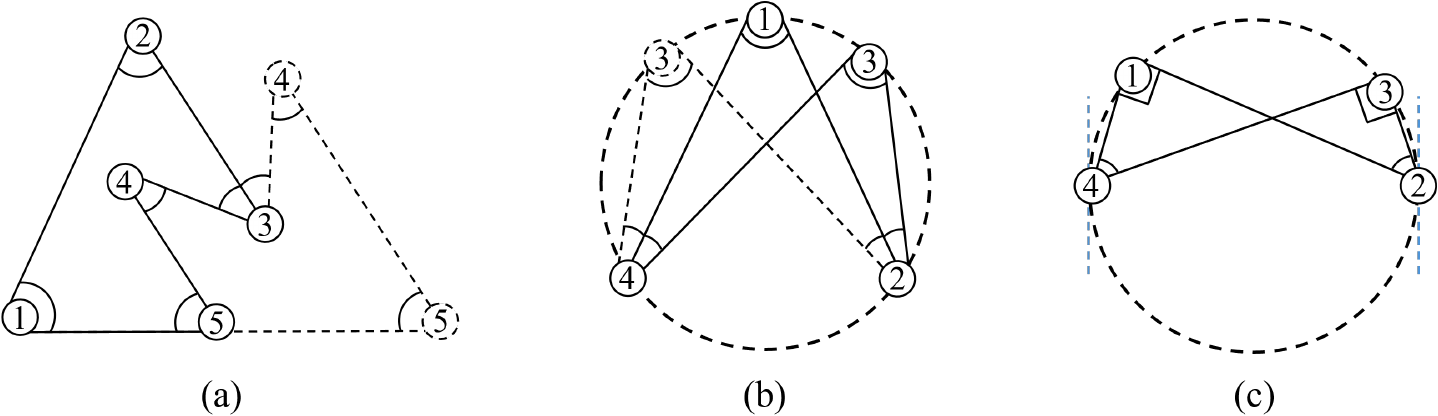} 
		\caption{(a) is a concave framework. (b) and (c) are non-planar frameworks and fall on two circles. In (a), dashed line $(4,5)$ and solid lines $(2,3)$ and $(4,5)$ are parallel, dashed line $(3,4)$ and solid line $(3,4)$ are symmetrical w.r.t. solid line $(2,3)$. In (b), dashed line $(2,3)$ (resp. dashed line $(3,4)$) and solid line $(2,3)$ (resp. solid line $(3,4)$) are symmetrical w.r.t. solid line $(1,2)$ (resp. solid line $(1,4)$). In (c), $\theta_{143}=\theta_{123}>\frac{\pi}{4}$ and $\theta_{214}=\theta_{234} = \frac{\pi}{2}$, and blue dashed lines indicate tangent lines at vertices 2 and 4, respectively.}
		\label{contradiction_for_5_cycle}
		\centering
	\end{figure}  

    \begin{lemma} \label{theorem: bearing and angle for C_n}
        Given a cycle framework $(\mathcal{C}_n,p)$ in $\mathbb{R}^2$, $ A^{-1}_{\mathcal{C}_n}(A_{\mathcal{C}_n}(p))$ is equivalent to $B^{-1}_{\mathcal{C}_n}(B_{\mathcal{C}_n}(p))$ up to uniform orthogonal transformations if $(\mathcal{C}_n,p)$ is convex.
    \end{lemma}

    \begin{lemma} \label{theorem: having same R for positively dependent}
    Given a cycle framework $(\mathcal{C}_4,p)$ in $\mathbb{R}^2$, $A^{-1}_{\mathcal{C}_4}(A_{\mathcal{C}_4}(p))$ is equivalent to $B^{-1}_{\mathcal{C}_4}(B_{\mathcal{C}_4}(p))$ up to uniform orthogonal transformations if $(\mathcal{C}_4,p)$ is planar.
    \end{lemma}

    The proofs are deferred to Appendix \ref{poof of theorem: bearing and angle for C_n} and Appendix \ref{poof of theorem: having same R for positively dependent}. Lemmas \ref{theorem: bearing and angle for C_n} and \ref{theorem: having same R for positively dependent} establish a nontrivial fact: angle measurements obtained in different local coordinate frames have the ability to collectively determine global bearings up to a uniform orthogonal transformation. A similar result for signed angles has been established in \cite[Lemma~3]{huang2026signed}, where the underlying graph is only required to be connected. However, the angles considered in this paper carry less information than signed angles, and hence may admit local angle-preserving reflections to prevent the conversion from angles to bearings, see Fig.~\ref{contradiction_for_5_cycle}(a) as an example. Lemmas \ref{theorem: bearing and angle for C_n} and \ref{theorem: having same R for positively dependent} exploit the cycle structure to exclude these undesired ambiguities.
    
    Some examples suggest that Lemmas \ref{theorem: bearing and angle for C_n} and \ref{theorem: having same R for positively dependent} can be extended to concave or non-planar cycle frameworks unless additional conditions on the angles are imposed. As shown in Figs.~\ref{contradiction_for_5_cycle}(a) and \ref{contradiction_for_5_cycle}(b), the set $A^{-1}_{\mathcal{C}_n}(A_{\mathcal{C}_n}(p))$ is not equivalent to $B^{-1}_{\mathcal{C}_n}(B_{\mathcal{C}_n}(p))$ up to uniform orthogonal transformations. However, the relationship holds when vertices 2, 3, 4, and 5 in Fig.~\ref{contradiction_for_5_cycle}(a) are collinear and the angles in Fig. \ref{contradiction_for_5_cycle}(b) satisfy $\theta_{143}=\theta_{123}>\frac{\pi}{4}$ and $\theta_{214}=\theta_{234} = \frac{\pi}{2}$ (i.e., Fig. \ref{contradiction_for_5_cycle}(c)). Next, we extend the relationship established for separate cycles to a successive union of cycles.

    The frameworks $(\bar{\mathcal{G}}_l, \bar{p}_l)$ for $l\in\{1,2\}$ satisfy $\bar{p}^i_1 = \bar{p}^i_2$ for $i\in \bar{\mathcal{V}}_1\bigcap\bar{\mathcal{V}}_2$, where $\bar{p}_l = [\dots,\bar{p}^i_l{}^\top,\dots]^\top$ and $\bar{p}^i_l\in\mathbb{R}^2$ is the coordinate assigned to vertex $i$ of $\bar{\mathcal{G}}_l = (\bar{\mathcal{V}}_l, \bar{\mathcal{E}}_l)$. The union of them is defined as $(\bar{\mathcal{G}}_1, \bar{p}_1)\bigcup(\bar{\mathcal{G}}_2, \bar{p}_2) \triangleq \left(\bar{\mathcal{G}}_1 \bigcup\bar{\mathcal{G}}_2, \bar{p}_1\bigcup\bar{p}_2\right)$ where $\bar{\mathcal{G}}_1 \bigcup\bar{\mathcal{G}}_2 \triangleq \left(\bar{\mathcal{V}}_1 \bigcup\bar{\mathcal{V}}_2, \bar{\mathcal{E}}_1\bigcup\bar{\mathcal{E}}_2\right)$ and $\bar{p}_1\bigcup\bar{p}_2$ is the stack of the coordinates of the vertices w.r.t. $(\bar{\mathcal{G}}_1, \bar{p}_1)$ and $(\bar{\mathcal{G}}_2, \bar{p}_2)$.

 	\begin{figure}[t]
		\centering
		\includegraphics[width=1\linewidth]{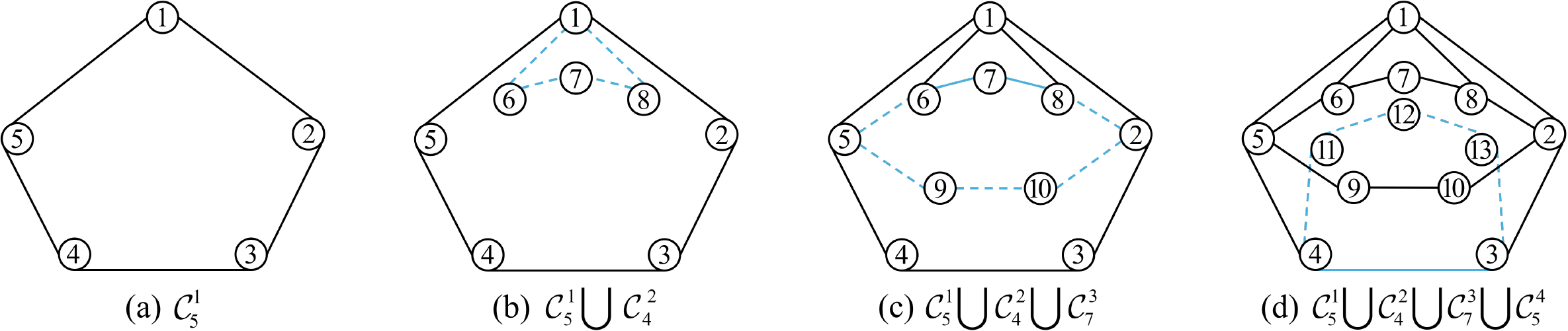} 
		\caption{Cycle compatible frameworks where blue lines (resp. dashed lines) indicate added cycles (resp. added edges). (a) $(\mathcal{C}^{1}_{5}, ^1\hspace{-1mm}\bar{p})$ is strictly convex. (b) $(\mathcal{C}^{2}_{4}, ^2\hspace{-0.5mm}\bar{p})$ is strictly concave. (c) $\mathcal{C}^{3}_{7}$ does not contain a path $\mathcal{P}$ with $|\mathcal{E}_P| \geq 4$ intersecting $\bar{\mathcal{G}}_{2}$ only at endpoints, and hence $b^\top_{59}(^3\hspace{-0.2mm}\bar{p})b_{2\hspace{0.5mm}10}(^3\hspace{-0.2mm}\bar{p}) \neq \pm 1$. (d) $(\mathcal{C}^{4}_{5}, ^4\hspace{-0.5mm}\bar{p})$ is required to be strictly convex since $\mathcal{C}^{4}_{5}$ contains a path with $|\mathcal{E}_P| \geq 4$ intersecting $\bar{\mathcal{G}}_{3}$ only at its endpoints 3, 4 and $|\mathcal{E}_P| = 4$.}
		\label{illustration_of_cycle_composition}
		\centering
	\end{figure}

    \begin{definition} \label{cycle fixable}
         A framework $(\bar{\mathcal{G}}_{c},\bar{p}_c)$ is said to be cycle compatible if it can be constructed recursively as follows. Start with a cycle framework $(\bar{\mathcal{G}}_{1},\bar{p}_1) = (\mathcal{C}^{1}_{n_{1}}, ^1\hspace{-1mm}\bar{p})$ which is strictly convex or strictly concave when $n_{1} = 4$, and strictly convex otherwise. $(\bar{\mathcal{G}}_{l},\bar{p}_{l}) = (\bar{\mathcal{G}}_{l-1},\bar{p}_{l-1}) \bigcup (\mathcal{C}^{l}_{n_{l}}, ^l\hspace{-1mm}\bar{p})$ with $l\geq2$ is formed by taking the union of $(\bar{\mathcal{G}}_{l-1},\bar{p}_{l-1})$ and the cycle framework $(\mathcal{C}^{l}_{n_{l}}, ^l\hspace{-1mm}\bar{p})$ subject to the following conditions:
        \begin{itemize}
        \item [(\rmnum{1})] $\mathcal{V}^{n_{l}}_C \bigcap \bar{\mathcal{V}}_{l-1} \neq \varnothing$ and $\bar{p}^i_{l-1} = ^l\hspace{-0.5mm}\hspace{-1mm}\bar{p}^i$ for $i\in \mathcal{V}^{n_{l}}_C \bigcap \bar{\mathcal{V}}_{l-1}$;
        \item [(\rmnum{2})] If $\mathcal{E}^{n_{l}}_C \bigcap \bar{\mathcal{E}}_{l-1} = \varnothing$, then $(\mathcal{C}^{l}_{n_{l}}, ^l\hspace{-1mm}\bar{p})$ is strictly convex or strictly concave when $n_{l} = 4$, and strictly convex otherwise;
        \item [(\rmnum{3})] If $\mathcal{E}^{n_{l}}_C \bigcap \bar{\mathcal{E}}_{l-1} \neq \varnothing$ and $\mathcal{C}^{l}_{n_{l}}$ contains a path $\mathcal{P}$ with $|\mathcal{E}_P| \geq 4$ intersecting $\bar{\mathcal{G}}_{l-1}$ only at endpoints, then $(\mathcal{C}^{l}_{n_{l}}, ^l\hspace{-1mm}\bar{p})$ is strictly convex;
        \item [(\rmnum{4})] If $\mathcal{E}^{n_{l}}_C \bigcap \bar{\mathcal{E}}_{l-1} \neq \varnothing$ and the path $\mathcal{P}$ in $\mathcal{C}^{l}_{n_{l}}$ intersecting $\bar{\mathcal{G}}_{l-1}$ only at endpoints satisfies $|\mathcal{E}_P| \leq 3$, $b^\top_{i_1i_2}(^l\hspace{-0.5mm}\bar{p}) b_{i_3i_4}(^l\hspace{-0.5mm}\bar{p}) \neq \pm 1$ for the path $\mathcal{P}$ with $\mathcal{E}_P \hspace{-1mm}=\hspace{-1mm} \{(i_1, i_2), (i_2,i_3), (i_3, i_4)\}$.
    \end{itemize}
    Here $\mathcal{C}^{l}_{n_{l}} = (\mathcal{V}^{n_{l}}_C,\mathcal{E}^{n_{l}}_C)$ for $l\in\{1,\dots,c\}$, $^l\hspace{-0.5mm}\bar{p} = [\dots,^l\hspace{-1mm}\bar{p}^i{}^\top,\dots]^\top$, and $^l\hspace{-0.5mm}\bar{p}^i\in\mathbb{R}^2$ is the coordinate mapped to the vertex $i\in \mathcal{V}^{n_{l}}_C$.
    \end{definition}

\begin{figure}
    \centering
    \includegraphics[width=0.8\linewidth]{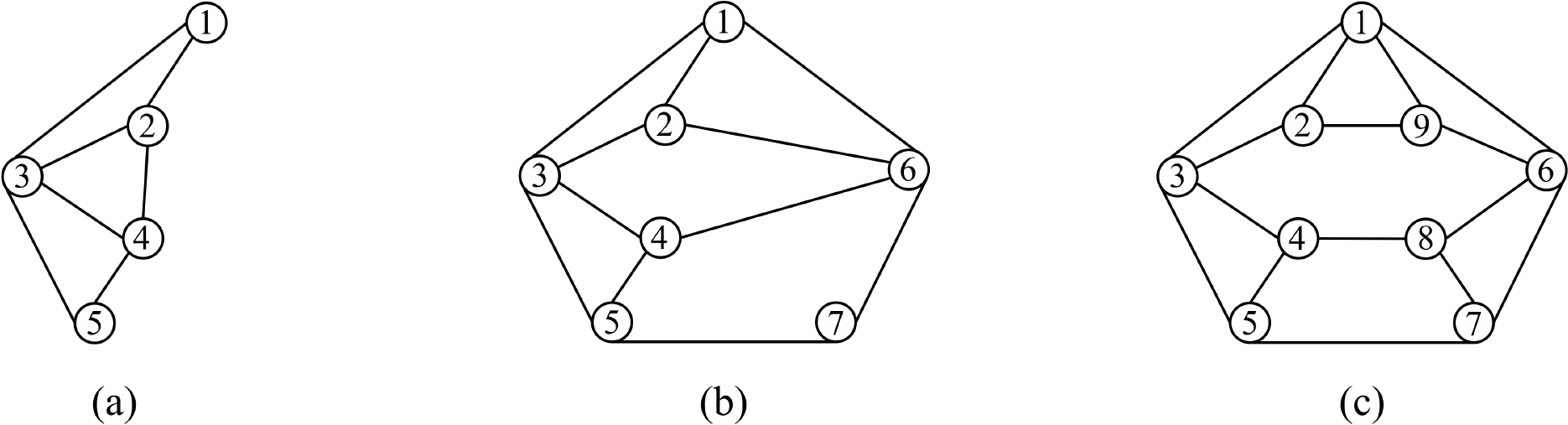}
    \caption{(a)-(c) are cycle compatible and IAR frameworks. (a) is a strongly non-degenerate triangulated framework constructed by vertex addition, where the selected vertices are required to be adjacent. (b) is a non-degenerate bilateration ordering constructed by vertex addition, where the selected vertices are not necessarily adjacent. The underlying graph of (c) is constructed using both vertex addition and edge splitting, and cannot be constructed using vertex addition alone.}
    \label{illustration_of_triangulation_bilateration_ordering}
    \centering
\end{figure}

    Fig. \ref{illustration_of_cycle_composition} provides examples of cycle compatible frameworks to show the geometric intuition behind Definition \ref{cycle fixable}. It should be noted that Definition \ref{cycle fixable} covers all possible numbers of shared vertices between $\bar{\mathcal{G}}_{l-1}$ and $\mathcal{C}^{l}_{n_{l}}$, i.e., $|\mathcal{V}^{n_{l}}_C\bigcap\bar{\mathcal{V}}_{l-1}| \in \{1, \cdots, n_l\}$. On the other hand, it can be known from Definition \ref{cycle fixable} that cycle compatible frameworks are more general than strongly non-degenerate triangulated frameworks \cite{jing2019angle} and non-degenerate bilateration orderings \cite{jing2021angle} which contain a hierarchical construction via vertex addition\footnote{Vertex addition and edge splitting are two parts of the Henneberg construction which provides a combinatorial characterization of Laman graphs $\mathcal{L}_n \triangleq (\mathcal{V}_n,\mathcal{E}_n)$ \cite{zhao2017laman}.}; see Fig.~\ref{illustration_of_triangulation_bilateration_ordering} for more illustrations.

\begin{lemma} \label{lemma: cycle fixable can determine the bearing}
    For any cycle compatible framework $(\mathcal{G},p)$, $q \in A^{-1}_{\mathcal{G}}(A_{\mathcal{G}}(p))$ if and only if $B_{\mathcal{G}}(q) = (I_{m}\otimes \mathscr{R})B_{\mathcal{G}}(p)$ for some $\mathscr{R}\in{\rm O}(2)$.
\end{lemma}

The proof is deferred to Appendix \ref{proof of lemma: cycle fixable can determine the bearing}. With the help of Lemma \ref{lemma: cycle fixable can determine the bearing}, we next show the equivalence between angle rigidity and bearing rigidity.

\begin{theorem} \label{theorem: equivalence in framework}
For any cycle compatible framework $(\mathcal{G},p)$, it is AR (GAR) if and only if it is BR (GBR).
\end{theorem}

The proof is deferred to Appendix \ref{proof of theorm: theorem: equivalence in framework}. Since a non-degenerate framework is GBR if and only if it is IBR \cite{huang2026signed} and cycle compatible frameworks must be non-degenerate, combining Theorem \ref{theorem: equivalence in framework} with \cite[Theorem 3]{zhao2015bearing} and \cite[Theorem 3]{jing2019angle}, we induce the following result.

\begin{theorem} \label{theorem: equivalence in framework2}
For any cycle compatible framework $(\mathcal{G},p)$, the rigidity properties AR, GAR, IAR, BR, GBR, and IBR are all equivalent. 
\end{theorem}

Theorem \ref{theorem: equivalence in framework2} indicates that a cycle compatible framework is angle fixable if and only if it is IAR, which can be examined by the rank of the angle rigidity matrix. 

\begin{theorem} \label{theorem: non-generic angle fixable framework}
    For any cycle compatible framework $(\mathcal{G},p)$, it is angle fixable if and only if it is IAR. 
\end{theorem}

It can be verified that the frameworks in Fig.~\ref{Sensor_network_illustration}(d) and Fig.~\ref{illustration_of_triangulation_bilateration_ordering} are both cycle compatible and IAR. Hence, by Theorem \ref{theorem: non-generic angle fixable framework}, the shape of these frameworks can be uniquely determined by angles, whereas the criteria proposed in \cite{jing2019angle, jing2021angle} are not applicable to the frameworks shown in Figs. \ref{Sensor_network_illustration}(d) and \ref{illustration_of_triangulation_bilateration_ordering}(c). However, determining whether a given framework is cycle compatible remains challenging. In the next subsection, we recognize angle fixable frameworks with generic configurations.

\begin{remark} \label{remark:10}
    The requirement for strictly convex (strictly concave) cycle frameworks in Definition \ref{cycle fixable} is not necessary. Let $\mathcal{P}$ be the longest path in $\mathcal{C}^{l}_{n_{l}}$ with $l\geq 2$ such that $\mathcal{P}$ intersects $\bar{\mathcal{G}}_{l-1}$ only at its endpoints. Convex or concave cycle frameworks are admissible for cycle compatibility if additional conditions are imposed: (\rmnum{1}) When $|\mathcal{E}_P|< 4$, each vertex $j\in\mathcal{V}^{n_{l}}_C \bigcap \bar{\mathcal{V}}_{l-1}$ has two neighbors $i,k$ in $(\bar{\mathcal{G}}_{l-1},\bar{p}_{l-1})$ such that $i, j, k$ are not collinear; (\rmnum{2}) When $|\mathcal{E}_P| \geq 4$ or $\mathcal{E}^{n_{l}}_C \bigcap \bar{\mathcal{E}}_{l-1} = \varnothing$, there is a vertex $j\in\mathcal{V}^{n_{l}}_C \bigcap \bar{\mathcal{V}}_{l-1}$ with two non-collinear edges in each of $(\mathcal{C}^{l}_{n_{l}},\bar{p}^l_{n_l})$ and $(\bar{\mathcal{G}}_{l-1},\bar{p}_{l-1})$, where the two pairs of edges are not identical.
\end{remark}

\subsection{Angle fixable graphs} \label{Angle fixability for graphs} Given a configuration $p\in\mathbb{R}^{2n}$, it is said to be \textit{generic} if the $2n$ elements of $p$ are algebraically independent \cite{jing2018weak}. The framework $(\mathcal{G}, p)$ is said to be a \textit{generic framework} if $p$ is generic. As shown in \cite{zhao2017laman} and \cite{jing2019angle}, $(\mathcal{G},p)$ is IAR for every generic configuration $p\in\mathbb{R}^{2n}$ if and only if $\mathcal{G}$ contains a Laman spanning subgraph. If $\mathcal{G}$ contains such a subgraph, since the set of generic configurations in $\mathbb{R}^{2n}$ has full measure, the configuration set of frameworks embedded by the graph $\mathcal{G}$ but not IAR is of measure zero. Under this fact, we say that IAR is a generic property of the graph, and a graph $\mathcal{G}$ is said to be IAR in $\mathbb{R}^2$ if a generic framework embedded by $\mathcal{G}$ is IAR in $\mathbb{R}^2$. 

In this subsection, we establish the equivalence between AR (GAR) and BR (GBR) for generic frameworks, based on which we show AR (GAR) is a generic property of the graph and obtain the equivalent graphical condition of angle fixability.

\begin{lemma} \label{lemma: angle can determine the bearing}
    Given a 2-edge-connected graph $\mathcal{G}$ and a generic configuration $p \in\mathbb{R}^{2n}$, $q \in A^{-1}_{\mathcal{G}}(A_{\mathcal{G}}(p))$ if and only if $B_{\mathcal{G}}(q) = (I_{m}\otimes \mathscr{R})B_{\mathcal{G}}(p)$ for some $\mathscr{R}\in{\rm O}(2)$.
\end{lemma}

The proof is deferred to Appendix \ref{poof of lemma: angle can determine the bearing}, in which it is shown that a local angle-preserving reflection in a cycle structure leads to a contradiction with the generic configuration. Note that Lemma \ref{lemma: angle can determine the bearing} is not necessarily valid for arbitrary configurations. For example, the framework in Fig.~\ref{contradiction_for_5_cycle}(a) has a cycle as its underlying graph, but it admits a local angle-preserving reflection to destroy the relationship shown in Lemma \ref{lemma: angle can determine the bearing}. Although the similar result in signed angle rigidity theory holds for any framework embedded by a connected graph \cite{huang2026signed}, signed angles inherently rely on a consistent orientation convention, which is not required for the angles considered here. On the other hand, a graph is 2-edge-connected if and only if it contains no bridge \cite{bickle2020fundamentals}. Therefore, any of the rigidity properties AR, GAR, or IAR for a framework $(\mathcal{G},p)$ requires the underlying graph $\mathcal{G}$ to be 2-edge-connected; otherwise, the two components separated by the bridge admit a nontrivial infinitesimal angle motion along the bridge. In particular, Laman graphs are 2-edge-connected.

Similar to Theorem \ref{theorem: equivalence in framework}, with the help of Lemma \ref{lemma: angle can determine the bearing}, we can establish the equivalence between AR (GAR) and BR (GBR) in the context of generic frameworks.

\begin{theorem}
    \label{theorm: angle rigidity is equivalent to bearing rigdity2}
    In $\mathbb{R}^2$, a generic framework $(\mathcal{G},p)$ is AR (GAR) if and only if it is BR (GBR).
\end{theorem}

Since a framework is either GBR or not GBR for all generic configurations \cite{huang2026signed}, it follows from Theorem \ref{theorm: angle rigidity is equivalent to bearing rigdity2} that AR (GAR) is a generic property of the graph. 

\begin{definition}
    In $\mathbb{R}^2$, the graph $\mathcal{G}$ is said to be angle fixable (AR, GAR, IAR) if there is a generic framework embedded by $\mathcal{G}$ is angle fixable (AR, GAR, IAR).
\end{definition}

Following \cite{jing2019angle}, we can introduce the definition of minimal angle fixability (AR, GAR, IAR) for graphs\footnote{A framework $(\mathcal{G},p)$ in $\mathbb{R}^2$ is said to be minimally AR (GAR, IAR) if there is no AR (GAR, IAR) spanning subframework $(\bar{\mathcal{G}},p)$ such that $|\bar{\mathcal{E}}|<|\mathcal{E}|$.}, accordingly. Combining Theorem \ref{theorm: angle rigidity is equivalent to bearing rigdity2} with \cite[Theorem 6]{huang2026signed} and \cite[Theorem 3]{jing2019angle}, we obtain the following result.
\begin{theorem} \label{Corollary: infinitesimal angle rigidity inplies global rigdity}
    In $\mathbb{R}^2$, a graph $\mathcal{G}$ is minimally angle fixable (AR, GAR, IAR) if and only if it is a Laman graph.
\end{theorem}

Theorem \ref{Corollary: infinitesimal angle rigidity inplies global rigdity} provides a necessary and sufficient graphical condition for determining the angle fixability of generic frameworks. Since the set of non-generic configurations in $\mathbb{R}^{2n}$ has measure zero, it is applicable to almost all configurations arising in various scenarios. It also provides an equivalent graphical condition beyond the existing sufficient criteria based on triangulation \cite{jing2019angle} or bilateration ordering \cite{jing2021angle}.

\begin{figure}[htbp]
		\centering
		\includegraphics[width=0.9\linewidth]{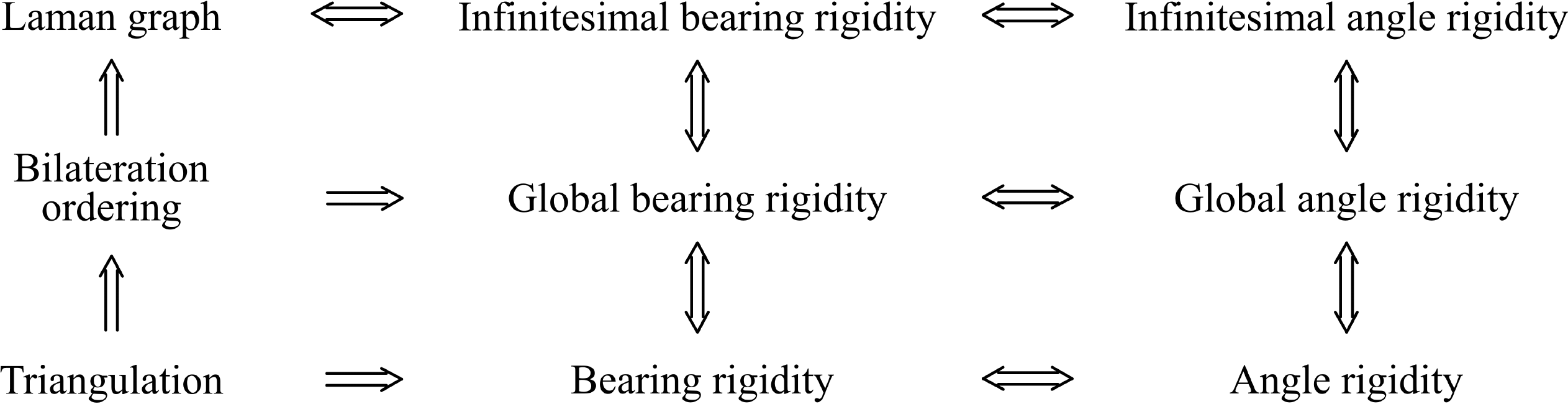}
		\caption{Relationship between bearing rigidity theory and angle rigidity theory for graphs in 2-D.}
		\label{relaiton conclusion graph2}
		\centering
	\end{figure}

Now, we have solved Problem \ref{problem: 1}. Fig. \ref{relaiton conclusion graph2} summarizes the relationships between different concepts of bearing rigidity and angle rigidity for graphs (or cycle compatible frameworks). These results provide the theoretical foundation for Problem \ref{problem: 2}.

\section{Application to sensor network localization}\label{sec:4}
In this section, by leveraging the theoretical results developed earlier, we first characterize the underlying graph of the sensor network to ensure the angle localizability. Building on this characterization, we further propose a distributed localization algorithm with guaranteed global convergence to solve Problem \ref{problem: 2}.

\subsection{Network characterization for distributed localization}To facilitate the design of a distributed localization algorithm, we propose Algorithm \ref{alg:1} to construct a class of planar Laman graphs, which will be exploited to determine the information flow and local update rules of the localization algorithm. Here, the notions of plane graphs, planar graphs, faces, face boundaries, and (weak) dual graphs are reviewed in Appendix \ref{Notation of graph theory}. 

\begin{algorithm}
    \renewcommand{\algorithmicrequire}{\textbf{Input:}}
    \renewcommand{\algorithmicensure}{\textbf{Output:}}
    \caption{Construction of a planar Laman graph $\mathcal{L}_h$}
    \label{alg:1}
    \begin{algorithmic}[1]
        \REQUIRE The number $h\geq3$ of vertices.
        \ENSURE $\mathcal{L}_h$
        \STATE Start with a plane graph $\tilde{\mathcal{G}}$ where $\mathcal{G}$ is a triangle with $\mathcal{V}=\{1,2,3\}$; let $l \leftarrow 1$.
        \FOR{$v = 4,\dots,h$}
            \STATE Perform one of the following two operations:
            \STATE \textbf{Vertex addition:} Choose two vertices $i,j \in \mathcal{V}$ that have a path of length at most $2$ along the boundary $\partial(\mathcal{F}_{\rm out})$ of $\tilde{\mathcal{G}}$. Adding vertex $v$ and two edges $(i,v),(j,v)$ to $\tilde{\mathcal{G}}$ creates a new face $\mathcal{F}^\prime$ with $|\partial(\mathcal{F}^\prime)| \leq 4$ and $\mathcal{F}^\prime\neq \mathcal{F}_{\rm out}$. 

            \STATE \textbf{Edge splitting:} Choose two faces $\mathcal{F}_1$ and $\mathcal{F}_2$ sharing an edge $(i,j) \in \mathcal{E}$ such that $|\partial(\mathcal{F}_1)| = 3$ in $\tilde{\mathcal{G}}$. Select a vertex $k$ such that there is a path from $i$ to $k$ of length at most $2$ along $\partial(\mathcal{F}_2)$. Remove edge $(i,j)$ in $\tilde{\mathcal{G}}$ to combine $\mathcal{F}_1$ and $\mathcal{F}_2$ into a single face. Adding vertex $v$ and three edges $(i,v),(j,v), (k,v)$ to $\tilde{\mathcal{G}}$ divides the merged face into three new faces $\mathcal{F}^\prime_i, i\in\{1,2,3\}$.
        \ENDFOR
        \STATE \textbf{return} The underlying graph $\mathcal{G}$ of the plane graph $\tilde{\mathcal{G}}$.
    \end{algorithmic}  
\end{algorithm}

    \begin{figure}
		\centering
		\includegraphics[width=1\linewidth]{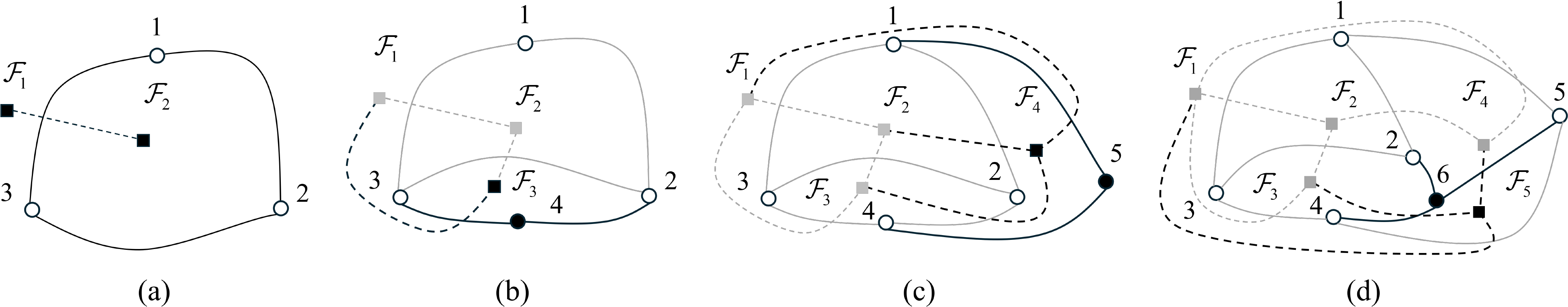}
		\caption{An illustration of Algorithm \ref{alg:1}, where black points and lines indicate the added elements. In (a)--(d), the graphs denoted by solid dots and lines are plane Laman graphs, and dual graphs are composed of square points for faces and dashed lines for capturing the adjacency of faces. (b)--(c) are formed by performing vertex addition. (d) is formed by performing edge splitting.}
		\label{illustration_of_algorithm_1}
		\centering
	\end{figure}

Fig. \ref{illustration_of_algorithm_1} illustrates the procedure of Algorithm \ref{alg:1}. It can be observed from the Lines 4--5 of Algorithm \ref{alg:1} that the introduction of new edges at each step divides a face and does not render the updated graph non-planar. Moreover, by Euler's formula for planar graphs and the fact that Laman graphs have $2n-3$ edges, any plane Laman graph must contain at least one triangle. Therefore, the two operations in Algorithm \ref{alg:1} can be performed at each iteration.

\begin{theorem} \label{theorem: proof for algorithm}
    The planar Laman graph $\mathcal{L}_h$ generated by Algorithm \ref{alg:1} can be obtained through a successive union of cycles $\mathcal{L}_h = \bigcup^c_{l=1}\mathcal{C}^{l}_{n_l}$, where $|\mathcal{E}^{n_1}_C| = 3$ and $3\leq |\mathcal{E}^{n_l}_C|\leq4$ and $\mathcal{E}^{n_{l}}_C \bigcap \bar{\mathcal{E}}_{l-1} \neq \varnothing$ for $\bar{\mathcal{E}}_{l-1} = \bigcup^{l-1}_{r=1}\mathcal{E}^{n_r}_{C} (l\geq2)$.
\end{theorem}

The proof is deferred to Appendix \ref{proof of theorem: proof for algorithm}, from which it can be known that the weak dual graph $\hat{\mathcal{G}} = (\hat{\mathcal{V}},\hat{\mathcal{E}})$ of the plane graph $\tilde{\mathcal{G}} $ produced by Algorithm \ref{alg:1} is connected and satisfies $|\partial(\mathcal{F}_i)| \leq 4$ for $\mathcal{F}_i\in \hat{\mathcal{V}}$. By extracting a spanning tree $(\hat{\mathcal{V}}_{\rm tree}, \hat{\mathcal{E}}_{\rm tree})$ from $\hat{\mathcal{G}}$, we can find a set of cycles $\mathcal{C}^l_{n_l}$ corresponding to elements $\mathcal{F}_l\in \hat{\mathcal{V}}_{\rm tree}$ and satisfying the conditions in Theorem \ref{theorem: proof for algorithm}.

Although the full set $\mathcal{T}_{\mathcal{G}}$ considered in Problem \ref{problem: 2} provides sufficient angle constraints, it is usually redundant. To reduce the angle measurement burden, we next exploit the graph produced by Algorithm \ref{alg:1} to construct a reduced set $\bar{\mathcal{T}}_{\mathcal{G}} \subset \mathcal{T}_{\mathcal{G}}$.

    \begin{figure}
		\centering
		\includegraphics[width=0.8\linewidth]{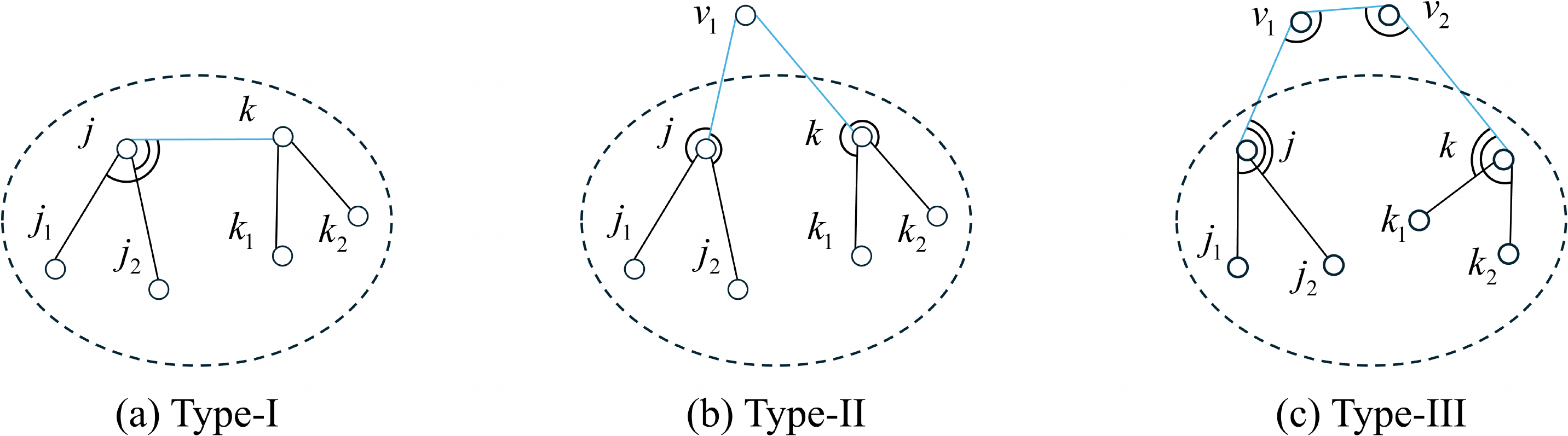}
		\caption{Illustration of introduced angle constraints. $v_i$ and the blue line segments denote the added vertex and edges, respectively. (a) The root edges of $(j,k)$ are $(j,j_1)$ and $(j,j_2)$. (b) The root edges of $(j,v_1)$ are $(j,j_1)$ and $(j,j_2)$. (c) The root edges of $(v_1,v_2)$ are $(v_1,j)$ and $(v_2,k)$.}
		\label{illustration_of_2_4_cycle_construction}
		\centering
	\end{figure}

\subsection{Construction of a reduced angle index set} Considering the cycles $\mathcal{C}^{1}_{n_1}, \dots, \mathcal{C}^{c}_{n_c}$ satisfying the conditions in Theorem \ref{theorem: proof for algorithm}, at the $l$-th union step (i.e., when $\mathcal{C}^{l}_{n_{l}}$ is introduced into the graph $\bigcup^{l-1}_{r=1}\mathcal{C}^{r}_{n_r}$ for $l\geq2$), we can introduce an angle constraint set denoted by $^l\bar{\mathcal{T}}_{\mathcal{G}} \subseteq \mathcal{T}_{\mathcal{G}}$ to connect the added edge with the root edges\footnote{At the $l$-th union step, the root edges of the added edge are the two edges selected to form angle constraints with it, which either both belong to $\bar{{\mathcal{E}}}_{l-1}$ or $\mathcal{E}^{n_l}_C\backslash\bar{{\mathcal{E}}}_{l-1}$.}; see Fig. \ref{illustration_of_2_4_cycle_construction} for more details. Then, the angle measurement set of the network is indexed by the set $\bar{\mathcal{T}}_{\mathcal{G}} =$$\mathcal{T}_{\mathcal{C}^1_{n_1}}\bigcup$$^2\bar{\mathcal{T}}_{\mathcal{G}}\bigcup\dots\bigcup$$^c\bar{\mathcal{T}}_{\mathcal{G}}$. 

Note that $\bar{\mathcal{T}}_{\mathcal{G}}$ may not be the minimal set for determining the geometric shape. Nevertheless, recalling the Henneberg construction, it can be verified by using induction that $|\mathcal{T}_{\mathcal{G}}| > 3+2[(2n-3)-3] = |\bar{\mathcal{T}}_{\mathcal{G}}|$ for $|\mathcal{V}|>3$, and hence, $\bar{\mathcal{T}}_{\mathcal{G}}$ is more efficient than the full set $\mathcal{T}_{\mathcal{G}}$. In what follows, we consider that each sensor $j\in\mathcal{V}$ has access to angle measurements associated with the set $\bar{\mathcal{T}}^j_{\mathcal{G}} = \{(i,j,k)\in\mathcal{V}^3:(i,j,k)\in\bar{\mathcal{T}}_{\mathcal{G}}\}$.

To ensure the solvability of the problem, we consider the following assumption: 

\begin{assumption} \label{ass: 1} The sensor network $(\mathcal{G},p,\mathcal{A})$ satisfies the following conditions: 
\begin{itemize}
    \item [(\rmnum{1})] $\mathcal{G}$ is formed by Algorithm \ref{alg:1} and $\{\mathcal{C}^{l}_{n_l}\}^c_{l=1}$ meets conditions in Theorem \ref{theorem: proof for algorithm};
    \item [(\rmnum{2})] The root edges of the edge $(i,j)\in\mathcal{E}\backslash\mathcal{E}^{n_1}_C$ are not parallel in $(\mathcal{G}, p)$;
    \item [(\rmnum{3})] $\mathcal{A} = \mathcal{V}^{n_1}_C$ and the anchors are not collinear in $(\mathcal{G}, p)$;
    \item [(\rmnum{4})] $(\mathcal{G}, p)$ is IAR.
\end{itemize}
\end{assumption}

It can be observed that the angle index graph\footnote{Each vertex of the angle index graph corresponds to an edge of $\mathcal{G}$ and two vertices are adjacent if the two corresponding edges are involved in an angle constraint belonging to $\bar{\mathcal{T}}_{\mathcal{G}}$ \cite{huang2026signed}.} associated with $\mathcal{G}$ and $\bar{\mathcal{T}}_{\mathcal{G}}$ is connected. According to the proof of \cite[Theorem 3]{jing2019angle}, the constraint set $\bar{\mathcal{T}}_{\mathcal{G}}$ is sufficient to support IAR. Moreover, it follows from Assumption \ref{ass: 1}(\rmnum{2}) and the proof of Lemma \ref{lemma: cycle fixable can determine the bearing} that for any $q\in \mathbb{R}^{2n}$ satisfying $b^\top_{ji}(q)b_{jk}(q) = b^\top_{ji}(p)b_{jk}(p)$, $(i,j,k)\in\bar{\mathcal{T}}_{\mathcal{G}}$, one has $b_{ji}(q) = \mathscr{R}b_{ji}(p)$, $(i,j)\in\mathcal{E}$ for some $\mathscr{R}\in{\rm O}(2)$, and hence, $b^\top_{ji}(q)b_{jk}(q) = b^\top_{ji}(p)b_{jk}(p)$ for $(i,j,k)\in\mathcal{T}_{\mathcal{G}}$. Consequently, $\bar{\mathcal{T}}_{\mathcal{G}}$ is able to act as $\mathcal{T}_{\mathcal{G}}$ for GAR. 

Under Assumption \ref{ass: 1}, it can be verified that $(\mathcal{G}, p)$ is cycle compatible. Since the anchors are not all collinear, it follows from Theorem \ref{theorem: non-generic angle fixable framework} that the network $(\mathcal{G},p,\mathcal{A})$ is angle localizable.

\subsection{Distributed localization algorithm} Motivated by \cite{zhao2016localizability}, each follower $j\in \mathcal{F}$ with an arbitrary condition $\hat{p}_j(0)\in\mathbb{R}^2$ updates its estimate according to
    \begin{eqnarray}
        \dot{\hat{p}}_j(t) &=& - \sum_{i\in\mathcal{N}_j}\left(I_2-(\|\hat{b}_{ji}(t)\|^2)^\dagger\hat{b}_{ji}(t)\hat{b}^\top_{ji}(t)\right)(\hat{p}_j(t)-\hat{p}_i(t))  \label{eq:location estimators1}.
    \end{eqnarray}
Here, $(\cdot)^\dagger$ denotes the Moore-Penrose pseudoinverse and is defined to satisfy $(0)^\dagger = 0$; $\hat{p}_i(t)$ is obtained through communication with neighbors $i\in\mathcal{N}_j$; $\hat{b}_{ji}(t) \in \mathbb{R}^{2}$ denotes the estimate of the global bearing $b_{ji}$, and is updated via the estimator (\ref{eq:location estimators2}) 

To estimate the bearing $b_{ji}$, sensor $j\in\mathcal{V}$ and its neighbor $i\in\mathcal{N}_j$ randomly generate a common initial estimate $\hat{b}_{ji}(0) = -\hat{b}_{ij}(0) \in \mathbb{R}^2$, and update $\hat{b}_{ji}(t)$ by
        \begin{eqnarray} 
            \hspace{-2mm}\dot{\hat{b}}_{ji}(t) &=& - {\rm sin}^2\hat{\beta}_{ji}\hat{b}_{ji}(t) + A_{ji}\hat{b}_{r_1k_1}(t) - B_{ji}\hat{b}_{r_2k_2}(t),\label{eq:location estimators2}
        \end{eqnarray}
        where $(r_1,k_1)$ and $(r_2,k_2)$ are the root edges of $(j,i)$, $\hat{\beta}_{ji}$ is the counterclockwise rotation angle from $\hat{b}_{r_1k_1}$ to $\hat{b}_{r_2k_2}$ with $\hat{\beta}_{ji} = 0$ whenever either of the root edge estimates is zero, and the coefficient matrices are
        \begin{eqnarray}
            A_{ji} &\triangleq& {\rm sin}\hat{\beta}_{ji}\mathscr{R}_o(\frac{\pi}{2}){\rm cos}\theta_{r_2k_2,ji} \in \mathbb{R}^{2\times 2}, \label{eq: estimator 4}\\
            B_{ji} &\triangleq& {\rm sin}\hat{\beta}_{ji}\mathscr{R}_o(\frac{\pi}{2}){\rm cos}\theta_{r_1k_1, ji} \in \mathbb{R}^{2\times 2}. \label{eq: estimator 5}
        \end{eqnarray}
        Here, ${\rm cos}\theta_{r_1k_1,ji} = b^\top_{r_1k_1}b_{ji}$, ${\rm cos}\theta_{r_2k_2,ji} = b^\top_{r_2k_2}b_{ji}$, $\hat{b}_{r_1k_1}$, and $\hat{b}_{r_2k_2}$ can be obtained by sensor $j$ or communicating with its neighbors. In the case of $(i,j) \in \mathcal{E}^{n_{1}}_C$, $b_{ij}$ is computed from the locations of anchors $i$ and $j$; $\hat{b}_{ij}(t) = b_{ij}$ and $\dot{\hat{b}}_{ij}(t) = 0$.

        \begin{theorem} \label{theorem:localization theorem}
            Under Assumption \ref{ass: 1}, by implementing the distributed protocol (\ref{eq:location estimators1})--(\ref{eq:location estimators2}),  $\hat{p}_j(t)$ asymptotically converges to $p_j$ under arbitrary initial estimates $\{\hat{p}_j(0)\}_{j\in\mathcal{F}}$ and $\{\hat{b}_{ji}(0)\}_{(j,i)\in\mathcal{E}}$.
        \end{theorem}

\begin{figure}
		\centering
		\includegraphics[width=0.90\linewidth]{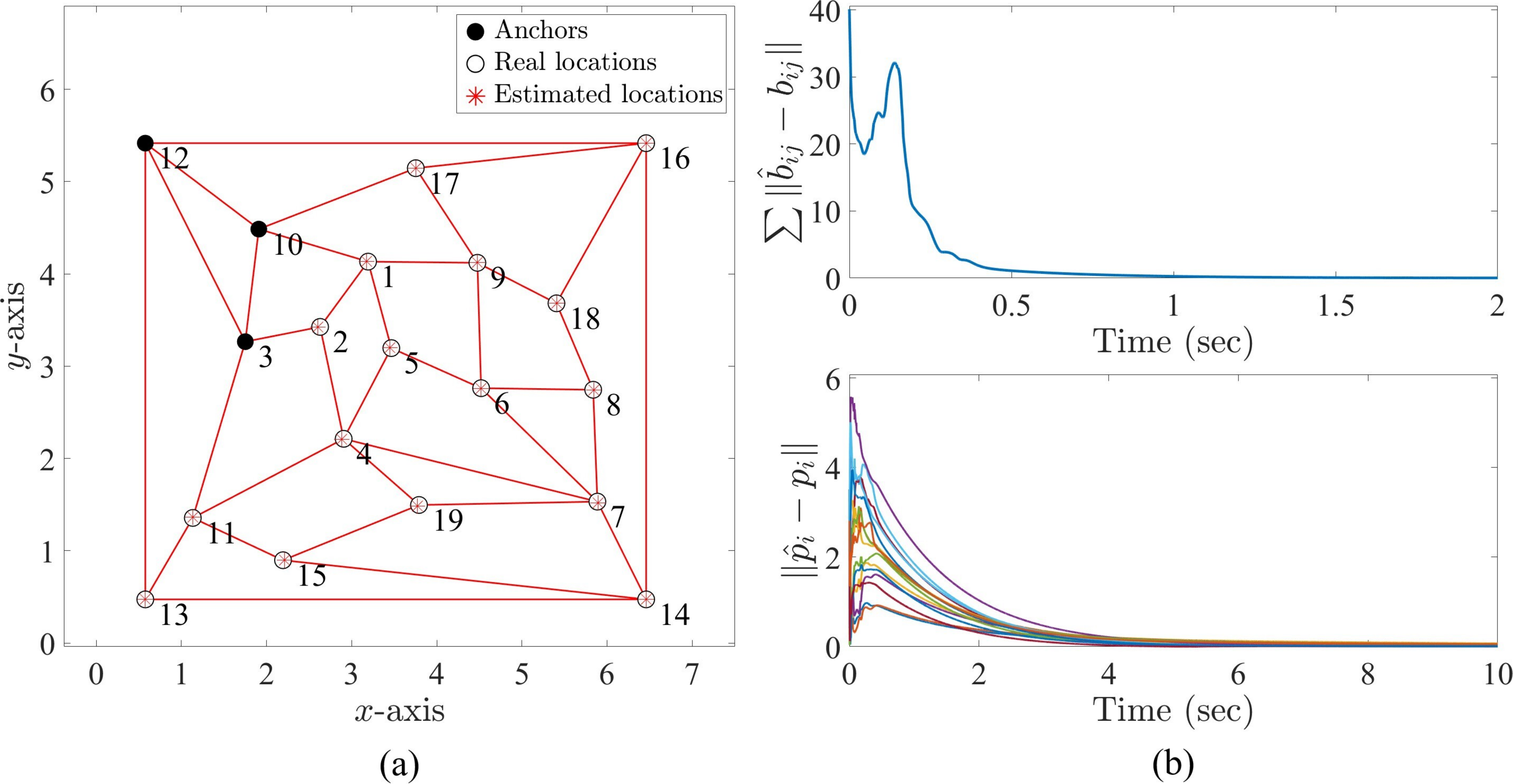}
		\caption{(a) The estimated locations converge to the real locations under the proposed protocol. (b) The evolutions of bearing and location estimation errors.}
		\label{unsigned_angle_localization}
		\centering
	\end{figure}

         The proof is deferred to Appendix \ref{proof of theorem:localization theorem}. It is worth noting that the sensing graphs $\mathcal{G}$ considered for the distributed protocol (\ref{eq:location estimators1})--(\ref{eq:location estimators2}) are more general than those in the literature. For example, the distributed protocol in \cite{jing2021angle} is not applicable to the network in Fig. \ref{unsigned_angle_localization}, but is admissible under the proposed distributed algorithm. Although the distributed protocol (\ref{eq:location estimators1})--(\ref{eq:location estimators2}) is analogous to the protocols in \cite{li2019globally, oh2026subtended}, which also perform network localization via an intermediate bearing-estimation stage, it accommodates different interpretations of the rotation directions across sensors. 
\begin{remark}
    The proposed protocol (\ref{eq:location estimators1})--(\ref{eq:location estimators2}) is applicable to any generic framework that contains a Laman spanning subgraph; however, the required bearing-estimation stage may incur additional communication and computational overhead. Consider the network $(\mathcal{G},p, \mathcal{A})$ with a generic $p$, $\mathcal{A} = \{1,6,7,8\}$, and $\mathcal{G}$ being the underlying graph in Fig. \ref{illustration_of_cycle_composition}(b). Since the bearings of edges $(1,5)$ and $(1,2)$ can be pre-determined from the angle measurements and the anchor locations, the edges $(1,5)$ and $(1,2)$ can be chosen as the root edges for $(4, 5)$. Additional communication is required for sensor $5$ to collect all angles indexed by $\mathcal{T}_{\mathcal{C}^1_{5}}$. According to Lemma \ref{lemma: angle can determine the bearing}, the angle between $b_{54}$ and $b_{12}$ can be determined by verifying which of the $2\times(5-2)$ vector combinations satisfies the angle constraints, where there are two possible bearings for the edge $(i,j)\in \{(2,3), (3,4), (4,5)\}$. The additional communicational and computational costs increase with the size of the cycle.
\end{remark}

\section{Simulation result}\label{sec:5} In this section, we present a numerical example to show the effectiveness of the proposed localization protocol.

Consider a sensor network $(\mathcal{L}_{19},p,\mathcal{A})$ in $\mathbb{R}^2$ as shown in Fig.~\ref{unsigned_angle_localization}, where the anchors are indexed by $\mathcal{A} = \{3,12,10\}$ and $(\mathcal{L}_{19},p,\mathcal{A})$ satisfies Assumption \ref{ass: 1}. $\mathcal{L}_{19} = \mathcal{C}^1_{n_1}\bigcup\cdots\bigcup\mathcal{C}^{17}_{n_{17}}$, where the vertex sets of the cycles are $\mathcal{V}^{n_1}_C = \{3,12,10\},
\mathcal{V}^{n_2}_C = \{3,10, 1, 2\},
\mathcal{V}^{n_3}_C = \{1,10,17,9\}, 
\mathcal{V}^{n_4}_C \hspace{-0.3mm}=\hspace{-0.3mm} \{10, 12, 16, 17\},
\mathcal{V}^{n_5}_C \hspace{-0.3mm}=\hspace{-0.3mm} \{12, 3, 11, 13\},
\mathcal{V}^{n_6}_C =\hspace{-0.3mm} \{11, 15, 14, 13\},
\mathcal{V}^{n_7}_C =\hspace{-0.3mm} \{12, 16, 14, 13\},
\mathcal{V}^{n_8}_C  = \{3, 2, 4, 11\},
\mathcal{V}^{n_9}_C = \{11, 4, 19, 15\}, 
\mathcal{V}^{n_{10}}_C  \hspace{-0.3mm}=\hspace{-0.3mm} \{15, 19, 7, 14 \}, 
\mathcal{V}^{n_{11}}_C \hspace{-0.3mm}=\hspace{-0.3mm} \{19, 4, 7\},
\mathcal{V}^{n_{12}}_C = \{2, 1, 5, 4\},
\mathcal{V}^{n_{13}}_C = \{1, 9, 6, 5\},
\mathcal{V}^{n_{14}}_C \hspace{-0.3mm}=\hspace{-0.3mm} \{6, 9, 18, \\ 8\},
\mathcal{V}^{n_{15}}_C \hspace{-0.3mm}=\hspace{-0.3mm} \{9, 17, 16, 18\},
\mathcal{V}^{n_{16}}_C \hspace{-0.3mm}=\hspace{-0.3mm} \{4,5,6, 7\},
\mathcal{V}^{n_{17}}_C = \{6, 8, 7\}$.

Given arbitrary $\{\hat{p}_j(0)\}_{j\in\mathcal{F}}$ and $\{\hat{b}_{ji}(0)\}_{(j,i)\in\mathcal{E}}$, we implement the distributed protocol (\ref{eq:location estimators1})--(\ref{eq:location estimators2}). The resulting convergence behavior and the evolutions of the estimation errors are shown in Fig.~\ref{unsigned_angle_localization}, which indicates that the estimated locations converge asymptotically to the real locations, in agreement with Theorem \ref{theorem:localization theorem}.

\section{Conclusions and future work}\label{sec:6}
In this paper, we provided the conditions for frameworks to be angle fixable with the help of cycle compatibility, and showed that the shape of Laman graphs can be uniquely determined by angles for almost all embeddings in the plane. Building on these results, we further designed a distributed localization algorithm for sensor networks with a more general class of underlying graphs than those considered in the existing literature. The proposed algorithm uses only angle measurements and achieves global convergence.

Future work includes: \rmnum{1}) finding the equivalent condition on the configuration for angle fixability of non-generic frameworks containing a Laman subgraph, \rmnum{2}) finding minimal angle constraint sets to support angle rigidity.

\appendix 

\section{Bearing rigidity theory}\label{appendix: Bearing Rigidity Theory}
In this section, we briefly review some notions and results from bearing rigidity theory \cite{zhao2015bearing}.

For a framework $(\mathcal{G},p)$ in $\mathbb{R}^2$, the bearing constraints can be stacked to form the \textit{bearing rigidity function}
    \begin{eqnarray} \label{equation: bearing rigidity function}
        B_{\mathcal{G}}(p) = \left[\cdots,b^\top_{ji},\cdots\right]^\top, (j,i) \in \mathcal{E}.
    \end{eqnarray}
    The framework is said to be \textit{bearing rigid} if there exists a neighborhood $U_p$ of $p$ such that $U_p \bigcap B^{-1}_{\mathcal{G}}(B_{\mathcal{G}}(p)) = U_p \bigcap B^{-1}_{\mathcal{K}}(B_{\mathcal{K}}(p))$, and to be \textit{globally bearing rigid} if $B^{-1}_{\mathcal{G}}(B_{\mathcal{G}}(p)) = B^{-1}_{\mathcal{K}}(B_{\mathcal{K}}(p))$, where $B^{-1}_{\mathcal{G}}(B_{\mathcal{G}}(p))$ is the configuration set of frameworks with the same bearings as $(\mathcal{G}, p)$. 

    An instantaneously continuous motion $\dot{p} \triangleq v$ is said to be an \textit{infinitesimal bearing motion} if it satisfies $\dot{B}_{\mathcal{G}}(p) = R_B(p)v = 0$, where $R_B(p) \triangleq \frac{\partial B_{\mathcal{G}}(p)}{\partial p}$ is the \textit{bearing rigidity matrix}. A framework is said to be \textit{infinitesimally bearing rigid} if all of its infinitesimal bearing motions consist of uniform translations and scalings.

    It has been reported in \cite{zhao2015bearing} that the geometric shape of a framework can be uniquely determined by bearing constraints (up to translations and scalings) if and only if it is infinitesimally bearing rigid. Furthermore, the authors in \cite{zhao2017laman} established the equivalent graphical condition for infinitesimal bearing rigidity. That is, the framework embedded by a Laman graph $\mathcal{L}_n \triangleq (\mathcal{V}_n,\mathcal{E}_n)$ satisfies minimal infinitesimal bearing rigidity for almost all configurations.

\section{Proof of Lemma \ref{theorem: bearing and angle for C_n}} \label{poof of theorem: bearing and angle for C_n} 

    \begin{proof} Suppose $(\mathcal{C}_n,p)$ is convex. Since $B^{-1}_{\mathcal{C}_n}(B_{\mathcal{C}_n}(p)) \subseteq A^{-1}_{\mathcal{C}_n}(A_{\mathcal{C}_n}(p))$, it suffices to show that for any $q\in A^{-1}_{\mathcal{C}_n}(A_{\mathcal{C}_n}(p))$, it always holds $(I_n \otimes\mathscr{R})^{-1}q\in B^{-1}_{\mathcal{C}_n}(B_{\mathcal{C}_n}(p))$ for some $\mathscr{R} \in {\rm O}(2)$. That is, $b_{i,i+1}(q) = \mathscr{R}b_{i,i+1}(p)$ for $i\in\mathcal{V}^n_C$ and some $\mathscr{R} \in {\rm O}(2)$. To prove this, we consider the following two cases.

        Case 1. $(\mathcal{C}_n,p)$ is degenerate. The proof is similar to that of \cite[Theorem 4]{jing2019angle}.
        
        Case 2. $(\mathcal{C}_n,p)$ is non-degenerate. It follows from $A_{\mathcal{C}_n}(q) = A_{\mathcal{C}_n}(p)$ that $(\mathcal{C}_n,q)$ is non-degenerate as well. Without loss of generality, for all $i\in \mathcal{V}^n_C$, we suppose $b^\top_{i+1,i+2}(p)b^\bot_{i,i+1}(p) \geq 0 $. Hence, the signed angle $\beta_{i+1}(p)$ from $b_{i,i+1}(p)$ to $b_{i+1,i+2}(p)$ satisfies $\beta_{i+1}(p) = {\rm arccos} (b^\top_{i+1,i+2}(p)b_{i,i+1}(p)) \in [0,\pi]$ and $\beta_{i+1}(p) + \theta_{i+1}(p) = \pi$, where $\theta_{i+1}(p) = {\rm arccos} (b^\top_{i+1,i+2}(p)b_{i+1,i}(p))$. From the convexity of $(\mathcal{C}_n,p)$, one has
        \begin{eqnarray} \label{signed angle equition1}
            \sum^{n-1}_{i=0}\theta_{i+1}(p) = (n-2)\pi\quad\text{and}\quad\sum^{n-1}_{i=0}\beta_{i+1}(p) = n\pi- \sum^{n-1}_{i=0}\theta_{i+1}(p) = 2\pi.
        \end{eqnarray}
        It follows from $A_{\mathcal{C}_n}(q) = A_{\mathcal{C}_n}(p)$ that $\beta_{i+1}(q) = \epsilon_{i+1}\beta_{i+1}(p)$, $i\in\mathcal{V}^n_C$ for some $\epsilon_{i+1}\in\{-1,1\}$. Note that $(\mathcal{C}_n,q)$ is embedded by a cycle. Therefore,
        \begin{eqnarray} \label{signed angle equitions}
            \sum^{n-1}_{i=0}\beta_{i+1}(q) = \sum^{n-1}_{i=0}\epsilon_{i+1}\beta_{i+1}(p) = 2k\pi
        \end{eqnarray}
        for some $k\in\mathbb{Z}$. Together with the fact that $\sum^{n-1}_{i=0}\beta_{i+1}(p) = 2\pi$ and $\beta_{i+1}(p)\geq0$, (\ref{signed angle equitions}) leads to the following three cases: 
        \begin{itemize}
            \item Case 2a. $\sum^{n-1}_{i=0}\beta_{i+1}(q) = 2\pi$, $\epsilon_{i+1} = 1$ for all $i\in \mathcal{V}^n_C$;
            \item Case 2b. $\sum^{n-1}_{i=0}\beta_{i+1}(q) = -2\pi$, $\epsilon_{i+1} = -1$ for all $i\in \mathcal{V}^n_C$;
            \item Case 2c. $\sum^{n-1}_{i=0}\beta_{i+1}(q) = 0$, $\epsilon_{i+1} = 1$ for some $i\in \mathcal{V}^n_C$, and $\epsilon_{i+1} = -1$ otherwise.
        \end{itemize} 
         Next, we show Cases 2a--2b lead to the conclusion (i.e., $b_{i,i+1}(q) = \mathscr{R}b_{i,i+1}(p)$ for $i\in\mathcal{V}^n_C$ and some $\mathscr{R} \in {\rm O}(2)$), and Case 2c conflicts with the fact that there exist scalars $d_{i,i+1}>0$ such that $\sum^n_{i=1}d_{i,i+1}b_{i,i+1}(q) = 0$.
        
        Case 2a. Similar to the proof of \cite[Lemma 3]{huang2026signed}, from $\beta_{i+1}(q) = \beta_{i+1}(p), i=0,\dots,n-1$, one has $b_{i,i+1}(q) = \mathscr{R}_o(\theta)b_{i,i+1}(p)$ for all $i\in\mathcal{V}^n_C$ and some $\theta\in[0,2\pi)$. 
        
        Case 2b. Note that there exists $\theta\in[0,2\pi)$ such that 
        \begin{eqnarray} \label{signed angle equitions 2}
            b_{12}(q) = \mathscr{R}_o(\theta)b_{12}(p) = \mathscr{H}_{b^\bot_{12}(q)}\mathscr{R}_o(\theta)b_{12}(p).
        \end{eqnarray}
        It follows from $\beta_{i+1}(q) = -\beta_{i+1}(p), i=0,\dots,n-1$ that for any $j\in\mathcal{V}^n_C$, one has 
        \begin{eqnarray} \label{signed angle equitions 3}
            \mathscr{R}^\top_o(\eta_{j,j+1}(p))b_{j,j+1}(p) = b_{12}(p) \quad\text{and}\quad \mathscr{R}^\top_o(-\eta_{j,j+1}(p))b_{j,j+1}(q) = b_{12}(q)
        \end{eqnarray}
        where $\eta_{j,j+1}(p) = \sum^{n-1}_{i=0}c_{i+1}\beta_{i+1}(p)$ for some $c_{i+1}\in\{0,1\}$. Combined (\ref{signed angle equitions 3}) with (\ref{signed angle equitions 2}), it follows from \cite[Lemma 8.2]{jing2019angle} that 
        \begin{eqnarray} 
            b_{j,j+1}(q) &=& \mathscr{R}_o(-\eta_{j,j+1}(p))\mathscr{H}_{b^\bot_{12}(q)}\mathscr{R}_o(\theta)\mathscr{R}^\top_o(\eta_{j,j+1}(p))b_{j,j+1}(p) \notag \\
                         &=& \mathscr{R}_o(2\alpha^\bot_{12})\mathscr{R}_o(-\eta_{j,j+1}(p)){\rm diag}(-1,1)\mathscr{R}^\top_o(\eta_{j,j+1}(p))\mathscr{R}_o(\theta)b_{j,j+1}(p) \notag \\
                         &=& \mathscr{H}_{b^\bot_{12}(q)}\mathscr{R}_o(\theta)b_{j,j+1}(p),
        \end{eqnarray}
        where $\mathscr{R}_o(-\eta_{j,j+1}(p)){\rm diag}(-1,1)\mathscr{R}^\top_o(\eta_{j,j+1}(p)) = {\rm diag}(-1,1)$ and $\alpha^\bot_{12}$ is the signed angle from $x$-axis to $b^\bot_{12}(q)$ in the counter-clockwise direction. As a result, there exists $\mathscr{R} = \mathscr{H}_{b^\bot_{12}(q)}\mathscr{R}_o(\theta)$ for some $\theta\in[0,2\pi)$ such that $b_{i,i+1}(q) = \mathscr{R}b_{i,i+1}(p), i\in\mathcal{V}^n_C$.

        Case 2c. $\sum^{n-1}_{i=0}\beta_{i+1}(q) = 0$ indicates that for the state sequence $\{x_i\}^{n}_{i=0}$ with $x_0 = \alpha_{n1}(q) \in[0,2\pi)$ and $x_{i+1} = x_i + \beta_{i+1}(q)$, the initial and final states are equivalent. Let $M$ and $m$ be the maximum and minimum values of $\{x_i\}^{n}_{i=0}$, respectively. From $\sum^{n-1}_{i=0}|\beta_{i+1}(q)| = 2\pi$ and $\alpha_{n1}(q)\leq|x_i| < 2\pi + \alpha_{n1}(q)$ for $i\in\{0,\dots,n\}$, one has 
        \begin{eqnarray} 
            2(M-m) \leq \sum^{n-1}_{i=0}|\beta_{i+1}(q)| = 2\pi \Rightarrow M-m \leq \pi.
        \end{eqnarray}
        Thus, the bearings $b_{i,i+1}(q), i\in\mathcal{V}^n_C$ satisfy $\alpha_{i,i+1}(q) = x_{i}$ and lie in the same half-plane of $\mathbb{R}^2$. Since $(\mathcal{C}_n,q)$ is not degenerate, they are not on the boundary. As a result, there are no strictly positive scalars $d_{i,i+1}>0$ such that $\sum^n_{i=1}d_{i,i+1}b_{i,i+1}(q) = 0$.
    \end{proof}

\section{Proof of Lemma \ref{theorem: having same R for positively dependent}}\label{poof of theorem: having same R for positively dependent} To simplify the discussion, we first provide the following definition and lemmas.
    
    \begin{definition} \label{definition: positively dependent}
        $(\mathcal{C}_n,p)$ in $\mathbb{R}^2$ is said to be \textit{positively dependent} if there exist strictly positive scalars $d_{i,i+1}$ such that $\sum^n_{i=1}d_{i,i+1}b_{i,i+1} = 0$, and is said to be \textit{strictly positively dependent} if it is positively dependent and $b_{i,i+1} \neq \pm b_{i+1,i+2}$ for each $i \in \mathcal{V}^n_C$.
    \end{definition}
    
    Next, we provide a useful lemma with regard to the planar framework $(\mathcal{C}_4,p)$.
    \begin{lemma} \label{lemma: conditions for planar framework}
        In $\mathbb{R}^2$, a framework $(\mathcal{C}_4,p)$ is a planar framework if and only if for any $q\in A^{-1}_{\mathcal{C}_4}(A_{\mathcal{C}_4}(p))$,  $(\mathcal{C}_4,q)$ is a planar framework.
    \end{lemma}
    \begin{proof}
        The sufficiency can be verified directly. To prove necessity, we first consider the cases with some $\theta_{ijk}$ being zero. In such cases, for any $q\in A^{-1}_{\mathcal{C}_4}(A_{\mathcal{C}_4}(p))$, $(\mathcal{C}_4,q)$ is either a degenerate quadrilateral (see Fig. \ref{illustration_of_Lemma_proof}(a)) or a degenerate framework in $\mathbb{R}^2$, and the result holds by definition.
        
        Without loss of generality, we suppose that no $\theta_{ijk}$ is zero and there exists $q\in A^{-1}_{\mathcal{C}_4}(A_{\mathcal{C}_4}(p))$ such that $(\mathcal{C}_4,q)$ is non‑planar. Since $\theta_{214} + \theta_{123} + \theta_{234} +\theta_{143} < 2\pi$, $(\mathcal{C}_4,p)$ is concave or non-planar. If $(\mathcal{C}_4,p)$ is concave,  without loss of generality, suppose that $\theta_{234} = \theta_{214}+\theta_{143}+\theta_{123}$. It can be observed from $(\mathcal{C}_4,q)$ that  $\theta_{214}+\theta_{123} = \theta_{234}+\theta_{143}$; then one has $\theta_{123} = 0$ and a contradiction arises. As a result, $(\mathcal{C}_4,p)$ is non-planar.
    \end{proof}

    The proof of Lemma \ref{theorem: having same R for positively dependent} proceeds in two steps. We first show the connection between angle constraints and bearing constraints in the context of strictly positively dependent frameworks, and then prove the main result.
    
    \begin{lemma} \label{lemma: having same R for strictly positively dependent}
    Given a strictly positively dependent and planar framework $(\mathcal{C}_4,p)$ in $\mathbb{R}^2$, $q\in A^{-1}_{\mathcal{C}_4}(A_{\mathcal{C}_4}(p))$ if and only if $(I_4\otimes\mathscr{R})^{-1}q \in B^{-1}_{\mathcal{C}_4}(B_{\mathcal{C}_4}(p))$ for some $\mathscr{R} \in {\rm O}(2)$.
    \end{lemma}
    
    \begin{proof}
        Since $B^{-1}_{\mathcal{C}_4}(B_{\mathcal{C}_4}(p)) \subseteq A^{-1}_{\mathcal{C}_4}(A_{\mathcal{C}_4}(p))$, it suffices to show that for any $q\in A^{-1}_{\mathcal{C}_4}(A_{\mathcal{C}_4}(p))$, it always holds $({\rm I}_4 \otimes\mathscr{R})^{-1}q\in B^{-1}_{\mathcal{C}_4}(B_{\mathcal{C}_4}(p))$ for some $\mathscr{R} \in {\rm O}(2)$. Note that $(\mathcal{C}_4,q)$ is strictly positively dependent since 
        \begin{eqnarray}
         \label{eq: proof 1}
            b^\top_{i,i+1}(q)b_{i+1,i+2}(q) = b^\top_{i,i+1}(p)b_{i+1,i+2}(p) \neq \pm 1, ~~~~i\in\mathcal{V}^4_C. 
        \end{eqnarray}
        Let $G_{214} = (b_{12},b_{14})$. Since $b^\top_{12}(p)b_{14}(p) = b^\top_{12}(q)b_{14}(q)$, one has $G_{214}^\top(p)G_{214}(p) = G_{214}^\top(q)G_{214}(q)$. Together with $\rank(G_{214}(p)) = 2$, it follows from \cite[Lemma 8.1]{jing2019angle} that $b_{12}(p) = \mathscr{R}_1b_{12}(q)$ and $b_{14}(p) = \mathscr{R}_1b_{14}(q)$ for some $\mathscr{R}_1 \in {\rm O}(2)$. Similarly, one has
        \begin{eqnarray}
            b_{21}(p) = \mathscr{R}_2b_{21}(q)\hspace{-3mm}&\quad,\quad&\hspace{-3mm}b_{23}(p) = \mathscr{R}_2b_{23}(q), \label{equation: 1}\\
            b_{41}(p) = \mathscr{R}_4b_{41}(q)\hspace{-3mm}&\quad,\quad&\hspace{-3mm}b_{43}(p) = \mathscr{R}_4b_{43}(q), \label{equation: 2}
        \end{eqnarray}
        for some $\mathscr{R}_2, \mathscr{R}_4 \in {\rm O}(2)$. Since $b_{41}(p) = \mathscr{R}_1b_{41}(q) = \mathscr{R}_4b_{41}(q)$, it follows from \cite[Lemma 8.3]{jing2019angle} that $\mathscr{R}_1 = \mathscr{R}_4$ or $\mathscr{R}_1\mathscr{H}_{b^\bot_{41}(q)} = \mathscr{R}_4$. Similarly, it follows from $\mathscr{R}_1b_{12}(q) = \mathscr{R}_2b_{12}(q)$ that $\mathscr{R}_1 = \mathscr{R}_2$ or $\mathscr{R}_1\mathscr{H}_{b^\bot_{12}(q)} = \mathscr{R}_2$. Therefore, there may exist four cases: 
        \begin{itemize}
            \item Case 1. $\mathscr{R}_1 = \mathscr{R}_2$ and $\mathscr{R}_1\mathscr{H}_{b^\bot_{41}(q)} = \mathscr{R}_4$;
            \item Case 2. $\mathscr{R}_1\mathscr{H}_{b^\bot_{12}(q)} = \mathscr{R}_2$ and $\mathscr{R}_1 = \mathscr{R}_4$;
            \item Case 3. $\mathscr{R}_1\mathscr{H}_{b^\bot_{12}(q)} = \mathscr{R}_2$ and $\mathscr{R}_1\mathscr{H}_{b^\bot_{41}(q)} = \mathscr{R}_4$;
            \item Case 4. $\mathscr{R}_1 = \mathscr{R}_2 = \mathscr{R}_4$.
        \end{itemize}
        
        \noindent It suffices to show that Cases 1--3 cause a contradiction, or lead to $\mathscr{R}_1 = \mathscr{R}_2 = \mathscr{R}_4$.
        
        Case 1. $\mathscr{R}_1 = \mathscr{R}_2$ and $\mathscr{R}_1\mathscr{H}_{b^\bot_{41}(q)} = \mathscr{R}_4$. 
        
        Suppose $b_{23}(q)$ is not parallel to $b_{41}(q)$. It follows from $b_{32}(p) = \mathscr{R}_1b_{32}(q)$ and $b_{34}(p) = \mathscr{R}_1\mathscr{H}_{b^\bot_{41}(q)}b_{34}(q)$ and $ b^\top_{34}(p)b_{32}(p) =  b^\top_{34}(q)b_{32}(q)$ that 
        \begin{eqnarray} \label{equation: case 1}
            b^\top_{34}(p)b_{32}(p) = b^\top_{34}(q)\mathscr{H}_{b^\bot_{41}(q)}b_{32}(q) = b^\top_{34}(q)b_{32}(q).
        \end{eqnarray}
        Since $b_{23}(q)$ is not parallel to $b_{41}(q)$, (\ref{equation: case 1}) implies that $b_{34}(q)$ is parallel to $b_{41}(q)$, and then $\|b^\top_{34}(p)b_{41}(p)\| = \|b^\top_{34}(q)b_{41}(q)\| = 1$. Hence, it can be concluded that $(\mathcal{C}_4,p)$ is not strictly positively dependent, and a contradiction arises.

        Suppose $b_{23}(q)$ is parallel to $b_{41}(q)$. Let $\tilde{q} = (I_4\otimes\mathscr{R}^\top_1)p$. From Lemma \ref{lemma: conditions for planar framework}, $(\mathcal{C}_4,\tilde{q})$ is a planar and strictly positively dependent framework satisfying $b_{ji}(\tilde{q}) = \mathscr{R}^\top_1b_{ji}(p),(j,i)\in\mathcal{E}^4_C$. Then, one has $b_{12}(\tilde{q}) = b_{12}(q), b_{23}(\tilde{q}) = b_{23}(q)$, $b_{41}(\tilde{q}) = b_{41}(q)$, and $b_{34}(\tilde{q}) = \mathscr{H}_{b^\bot_{41}(q)}b_{34}(q)$.

    \begin{figure}
		\centering
		\includegraphics[width=.95\linewidth]{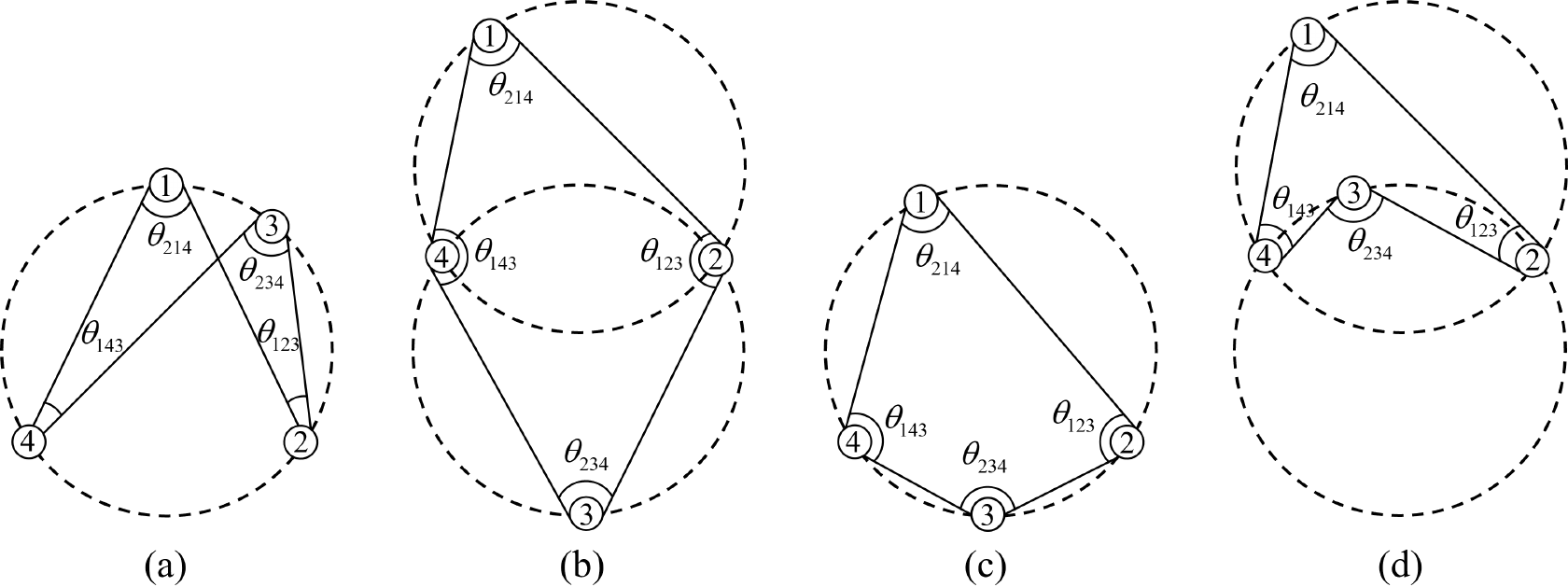} 
		\caption{The geometric interpretations of $\theta_{214} = \theta_{234}$ and $\theta_{214}+\theta_{234} = \pi$ are (a)--(b) and (c)--(d), respectively. (a) $(\mathcal{C}_4,q)$ is a non-planar framework and falls on a circle. (b) $(\mathcal{C}_4,q)$ falls on two congruent intersecting circles, in which vertex 2 and vertex 4 are at the intersection points. (c) $(\mathcal{C}_4,q)$ falls on a circle and is convex. (d) $(\mathcal{C}_4,q)$ falls on two congruent intersecting circles.}
		\label{illustration_of_Lemma_proof_3}
	\end{figure}        
    
        Since $\mathscr{H}_{b^\bot_{41}(q)}\mathscr{R}_o(\frac{\pi}{2}) = \mathscr{R}_o(\frac{\pi}{2})\mathscr{H}_{b_{41}(q)}$ and $\mathscr{H}_{b_{41}(q)}b_{41}(q) = -b_{41}(q)$, one has
        \begin{eqnarray} \label{equation:jjj}
            b^\top_{34}(\tilde{q})b^\bot_{41}(\tilde{q}) = b^\top_{34}(q)\mathscr{H}_{b^\bot_{41}(q)}\mathscr{R}_o(\frac{\pi}{2})b_{41}(q) = -b^\top_{34}(q)b^\bot_{41}(q).
        \end{eqnarray}
        Note that $(\mathcal{C}_4,q)$ is strictly positively dependent and $b^\top_{23}(q)b^\bot_{41}(q) = 0$. According to \cite[Lemma 1]{trinh2018formations}, $b^\top_{12}(q)b^\bot_{41}(q)$ and $b^\top_{34}(q)b^\bot_{41}(q)$ have opposite sign. Thus, $b^\top_{12}(\tilde{q})b^\bot_{41}(\tilde{q})\hspace{-0.2mm}=\hspace{-0.2mm} b^\top_{12}(p)b^\bot_{41}(p)$ and $b^\top_{34}(\tilde{q})b^\bot_{41}(\tilde{q}) \hspace{-0.3mm}=\hspace{-0.2mm} -b^\top_{34}(q)b^\bot_{41}(q)$ share the same sign. For any $i\in\mathcal{V}^4_C$, it holds that $b^\top_{i,i+1}(\tilde{q})b^\bot_{41}(\tilde{q}) \geq 0$ or $b^\top_{i,i+1}(\tilde{q})b^\bot_{41}(\tilde{q}) \leq 0$. Since $\tilde{q} = (I_4\otimes\mathscr{R}^\top_1)p$, according to \cite[Lemma 1]{trinh2018formations}, $(\mathcal{C}_4,\tilde{q})$ and $(\mathcal{C}_4,p)$ are not strictly positively dependent.

        Case 2. $\mathscr{R}_1\mathscr{H}_{b^\bot_{12}(q)} = \mathscr{R}_2$ and $\mathscr{R}_1 = \mathscr{R}_4$. The proof is similar to Case 1.

        Case 3. $\mathscr{R}_1\mathscr{H}_{b^\bot_{12}(q)} = \mathscr{R}_2$ and $\mathscr{R}_1\mathscr{H}_{b^\bot_{41}(q)} = \mathscr{R}_4$. 

        Let $\tilde{q} = (I_4\otimes\mathscr{R}^\top_1)p$. Then, $b_{12}(\tilde{q}) = b_{12}(q), b_{23}(\tilde{q}) = \mathscr{H}_{b^\bot_{12}(q)}b_{23}(q), b_{34}(\tilde{q}) = \mathscr{H}_{b^\bot_{41}(q)}b_{34}(q)$, and $b_{41}(\tilde{q}) = b_{41}(q)$. Following from \cite[Lemma 8.2]{jing2019angle}, $\mathscr{H}_{b^\bot_{41}(q)} = \mathscr{H}_{b^\bot_{14}(q)}$, and $\tilde{q}\in A^{-1}_{\mathcal{C}_4}(A_{\mathcal{C}_4}(q))$, one has
        \begin{eqnarray} \label{equation: 11}
            b^\top_{32}(\tilde{q})b_{34}(\tilde{q}) &=&  b^\top_{32}(q)\mathscr{H}_{b^\bot_{12}(q)}\mathscr{H}_{b^\bot_{41}(q)}b_{34}(q) \notag \\
            &=& b^\top_{32}(q)\mathscr{R}_o(2(\alpha_{12}(q)-\alpha_{14}(q)))b_{34}(q) \notag \\
            &=& b^\top_{32}(q)b_{34}(q).
        \end{eqnarray}
        Without loss of generality, we assume that vertices 1, 2, and 4 in $(\mathcal{C}_4,q)$ are labeled as in Fig. \ref{illustration_of_Lemma_proof_3}. It can be verified that $\mathscr{R}_o(2\theta_{214}) = \mathscr{R}_o(2(\alpha_{12}(q) - \alpha_{14}(q)))$, where $\theta_{214} = \arccos (b^\top_{12}(q)b_{14}(q))$. Together with (\ref{equation: 11}) and the fact that $\theta_{214} \neq 0,\pi$, it holds that $b^\top_{32}(q)\mathscr{R}_o(2\theta_{214})b_{34}(q) = b^\top_{32}(q)b_{34}(q)$ if and only if $\theta_{214} = \theta_{234}$ or $\theta_{214}+\theta_{234} = \pi$. Therefore, it suffices to show that $(\mathcal{C}_4,\tilde{q})$ is not planar and strictly positively dependent or $\mathscr{R}_1 = \mathscr{R}_2 = \mathscr{R}_4$. We consider the following cases.
        
    Case 3a. $\theta_{214} = \theta_{234}$ and $(\mathcal{C}_4,q)$ falls on a circle. 
        
        It can be observed from Fig. \ref{illustration_of_Lemma_proof_3}(a) that $(\mathcal{C}_4,q)$ is a non-planar framework. Hence, it follows from Lemma \ref{lemma: conditions for planar framework} and $\tilde{q}\in A^{-1}_{\mathcal{C}_4}(A_{\mathcal{C}_4}(q))$ that $(\mathcal{C}_4,\tilde{q})$ is not planar. 
        
    Case 3b. $\theta_{214} = \theta_{234}$ and $(\mathcal{C}_4,q)$ falls on two congruent and intersecting circles.
        
        Since $b^\top_{32}(q)\mathscr{R}_o(2\theta_{234})b_{34}(q) =  b^\top_{32}(q)b_{34}(q)$, $\theta_{234} = \frac{\pi}{2}$ or $\pi$ from Fig. \ref{illustration_of_Lemma_proof_3}(b). If $\theta_{234} = \pi$, it conflicts with the strictly positive dependency of $(\mathcal{C}_4,\tilde{q})$ directly. Thus, $\theta_{214} = \theta_{234} = \frac{\pi}{2}$, and then, the two circles passing through the vertices in $(\mathcal{C}_4,q)$ share a common center; vertices $1$ and $3$ lie on different arcs of the same circle. Hence, the cycle framework $(\mathcal{C}_4,q)$ is convex. From Lemma \ref{theorem: bearing and angle for C_n}, $\mathscr{R}_1 = \mathscr{R}_2 = \mathscr{R}_4$. 
        
    Case 3c. $\theta_{214}+\theta_{234} = \pi$ and $(\mathcal{C}_4,q)$ falls on a circle.
        
        Since the cycle framework $(\mathcal{C}_4,q)$ is convex, one has $\mathscr{R}_1 = \mathscr{R}_2 = \mathscr{R}_4$. 
    
    Case 3d. $\theta_{214}+\theta_{234} = \pi$ and $(\mathcal{C}_4,q)$ falls on two congruent intersecting circles.
    
    From (\ref{equation: 11}) and $\tilde{q}\in A^{-1}_{\mathcal{C}_4}(A_{\mathcal{C}_4}(q))$, one has
    \begin{eqnarray}
        b^\top_{32}(\tilde{q})b_{34}(\tilde{q}) &=& b^\top_{32}(q)\mathscr{R}_o(2(\pi-\theta_{234}))b_{34}(q) = b^\top_{32}(q)b_{34}(q),
    \end{eqnarray}
    which implies $\theta_{214} = \theta_{234} = \frac{\pi}{2}$ since $\theta_{234} \notin\{0,\pi\}$. Thus, the circles passing through the vertices of $(\mathcal{C}_4,q)$ share a common center, and vertices 1, 3 lie on the same arc of a circle; then, $(\mathcal{C}_4,q)$ is not planar. Hence, $(\mathcal{C}_4,\tilde{q})$ is not planar.
    \end{proof}

    With the help of Lemma \ref{lemma: having same R for strictly positively dependent}, we next prove Lemma \ref{theorem: having same R for positively dependent}.
\begin{proof}
    From Lemma \ref{lemma: having same R for strictly positively dependent}, the result holds when $(\mathcal{C}_4,p)$ is planar and strictly positively dependent. It suffices to show that if $b_{i,i+1}(p) = \pm b_{i+1,i+2}(p)$ for some $i\in\mathcal{V}^4_C$ , then for any $q\in A^{-1}_{\mathcal{C}_4}(A_{\mathcal{C}_4}(p))$, it holds that $b_{j,j+1}(q) = \mathscr{R}b_{j,j+1}(p)$ for $j\in\mathcal{V}^4_C$ and some $\mathscr{R} \in {\rm O}(2)$. To prove the statement, we consider the following two cases.

        Case 1. There is more than one $i\in\mathcal{V}_{C}$ with $b_{i,i+1}(p) = \pm b_{i+1,i+2}(p)$. Then, $(\mathcal{C}_4,p)$ is a degenerate framework. The proof is similar to that of \cite[Theorem 4]{jing2019angle}.

        Case 2. There exists only one $i\in\mathcal{V}_{C}$ with $b_{i,i+1}(p) = \pm b_{i+1,i+2}(p)$. Without loss of generality, let $b_{41}(p) = -b_{12}(p)$. Since $A_{\mathcal{C}_4}(p) = A_{\mathcal{C}_4}(q)$, it follows from $b^\top_{14}(p)b_{12}(p) = b^\top_{14}(q)b_{12}(q) = 1$ that $b_{14}(q) = b_{12}(q)$. Then, one has
        \begin{eqnarray}
            b^\top_{12}(p)b_{34}(p) = b^\top_{14}(p)b_{34}(p) = b^\top_{14}(q)b_{34}(q) = b^\top_{12}(q)b_{34}(q).
        \end{eqnarray}
        Let $G_{1234} = (b_{12},b_{23},b_{34})\in \mathbb{R}^{3\times 2}$. One has $G_{1234}^\top(p)G_{1234}(p) = G_{1234}^\top(q)G_{1234}(q)$. Note that $b^\top_{23}(p)b_{34}(p) = b^\top_{23}(q)b_{34}(q) \neq \pm1$. It follows from \cite[Lemma 8.1]{jing2019angle} that $b_{12}(p) = \mathscr{R}b_{12}(q)$, $b_{23}(p) = \mathscr{R}b_{23}(q)$, and $b_{34}(p) = \mathscr{R}b_{34}(q)$ for some $\mathscr{R} \in {\rm O}(2)$. Since $b_{14}(p) = b_{12}(p)$ and $b_{14}(q) = b_{12}(q)$, it follows that $b_{14}(p) = \mathscr{R}b_{14}(q)$.
\end{proof}

\section{Proof of Lemma \ref{lemma: cycle fixable can determine the bearing}} \label{proof of lemma: cycle fixable can determine the bearing}
\begin{proof}
        Sufficiency follows from $B^{-1}_{\mathcal{G}}(B_{\mathcal{G}}(p)) \subseteq A^{-1}_{\mathcal{G}}(A_{\mathcal{G}}(p))$. We next prove the necessity by induction. 
        
        For $c= 1$, the result holds directly from Lemmas \ref{theorem: bearing and angle for C_n} and \ref{theorem: having same R for positively dependent}.

        For $c \geq l > 1$, let $\bar{q}_{l} \in A^{-1}_{\bar{\mathcal{G}}_{l}}(A_{\bar{\mathcal{G}}_{l}}(\bar{p}_{l}))$. Suppose $b_{ji}(\bar{q}_{l}) = \mathscr{R}_1b_{ji}(\bar{p}_{l})$, $(j,i)\in \bar{\mathcal{E}}_{l-1}$ for some $\mathscr{R}_1\in {\rm O}(2)$. Let $\mathcal{P} = (\mathcal{V}_P,\mathcal{E}_P)$ be the longest path in $\mathcal{C}^{l}_{n_{l}}$ such that $\mathcal{P}$ intersects $\bar{\mathcal{G}}_{l-1}$ only at its endpoints. It suffices to show that $b_{ji}(\bar{q}_{l}) = \mathscr{R}_1b_{ji}(\bar{p}_{l})$ for $(j,i)\in \bar{\mathcal{E}}_{l}\backslash \bar{\mathcal{E}}_{l-1}$ under the following cases:

        Case 1. $|\mathcal{V}^{n_{l}}_C \bigcap \bar{\mathcal{V}}_{l-1}| \geq 2$ and $|\mathcal{E}_P| = 1$. Let $(j,i) \in \bar{\mathcal{E}}_{l} \backslash \bar{\mathcal{E}}_{l-1}$. Since each added cycle framework $(\mathcal{C}^i_{n_i},^i\hspace{-1mm}\bar{p})$ is strictly convex or strictly concave for $l\geq i\geq1$, there exist $k_1, k_2 \in \mathcal{N}^j_{\bar{\mathcal{G}}_{l-1}}$ such that $b_{jk_1}(\bar{p}_{l})$ is not parallel to $b_{jk_2}(\bar{p}_{l})$. Therefore, it follows from $(i,j,k_1), (i,j,k_2) \in \mathcal{T}_{\bar{\mathcal{G}}_{l}}$ that $b_{ji}(\bar{p}_{l})$ is uniquely determined by 
        \begin{eqnarray} \label{equiation: b_ji}
            b_{ji}(\bar{p}_{l}) &=& \begin{bmatrix}
            b^\top_{jk_1}(\bar{p}_{l}) \\
            b^\top_{jk_2}(\bar{p}_{l})
            \end{bmatrix}^{-1}\begin{bmatrix}
            {\rm cos} \theta_{ijk_1} \\
            {\rm cos} \theta_{ijk_2}
            \end{bmatrix}.
        \end{eqnarray}
        From $b_{jk_1}(\bar{q}_{l}) = \mathscr{R}_1b_{jk_1}(\bar{p}_{l})$ and $b_{jk_2}(\bar{q}_{l}) = \mathscr{R}_1b_{jk_2}(\bar{p}_{l})$, it can be known that $b_{jk_1}(\bar{q}_{l})$ is not parallel to $b_{jk_2}(\bar{q}_{l})$ as well, and $b_{ji}(\bar{q}_{l})$ satisfies
                \begin{eqnarray}
            b_{ji}(\bar{q}_{l}) &=& \begin{bmatrix}
            b^\top_{jk_1}(\bar{q}_{l}) \\
            b^\top_{jk_2}(\bar{q}_{l})
            \end{bmatrix}^{-1}\begin{bmatrix}
            {\rm cos} \theta_{ijk_1} \\
            {\rm cos} \theta_{ijk_2}
            \end{bmatrix} = \left\{\begin{bmatrix}
            b^\top_{jk_1}(\bar{p}_{l}) \\
            b^\top_{jk_2}(\bar{p}_{l})
            \end{bmatrix}\cdot \mathscr{R}^\top_1\right\}^{-1}\begin{bmatrix}
            {\rm cos} \theta_{ijk_1} \\
            {\rm cos} \theta_{ijk_2}
            \end{bmatrix}.
        \end{eqnarray}
    Combined with (\ref{equiation: b_ji}), one has $b_{ji}(\bar{q}_{l}) = \mathscr{R}_1b_{ji}(\bar{p}_{l})$ for any $(j,i) \in \bar{\mathcal{E}}_{l} \backslash \bar{\mathcal{E}}_{l-1}$.

    Case 2. $|\mathcal{V}^{n_{l}}_C \bigcap \bar{\mathcal{V}}_{l-1}| \geq 2$ and $|\mathcal{E}_P| = 2$. The proof is similar to Case 1.

    Case 3. $|\mathcal{V}^{n_{l}}_C \bigcap \bar{\mathcal{V}}_{l-1}| \geq 2$ and $|\mathcal{E}_P| = 3$. Without loss of generality, let $\mathcal{E}_P = \{(i_1,i_2), (i_2,i_3), (i_3,i_4)\}$. According to Case 1, one has $b_{i_2i_1}(\bar{q}_{l}) = \mathscr{R}_1b_{i_2i_1}(\bar{p}_{l})$ and $b_{i_3i_4}(\bar{q}_{l}) = \mathscr{R}_1b_{i_3i_4}(\bar{p}_{l})$. From Definition \ref{cycle fixable}, $b_{i_2i_1}(\bar{p}_{l})$ is not parallel to $b_{i_3i_4}(\bar{p}_{l})$. Similar to Case 1, it follows from $(i_1,i_2,i_3), (i_2,i_3,i_4) \in \mathcal{T}_{\bar{\mathcal{G}}_{l}}$ that $b_{i_2i_3}(\bar{q}_{l}) = \mathscr{R}_1b_{i_2i_3}(\bar{p}_{l})$. Together with the discussion in Cases 1--2, one has $b_{ji}(\bar{q}_{l}) = \mathscr{R}_1b_{ji}(\bar{p}_{l}), (j,i)\in \bar{\mathcal{E}}_{l}\backslash \bar{\mathcal{E}}_{l-1}$.
    
    Case 4. $|\mathcal{V}^{n_{l}}_C \bigcap \bar{\mathcal{V}}_{l-1}| \geq 2$ and $|\mathcal{E}_P| \geq 4$. We consider two cases.

    Case 4a. Suppose that there is no vertex $j \in \mathcal{V}^{n_{l}}_C \bigcap \bar{\mathcal{V}}_{l-1}$ having $(j,i) \in \mathcal{E}^{n_{l}}_C \bigcap \bar{\mathcal{E}}_{l-1}$ and $(j,k) \in \mathcal{E}^{n_{l}}_C \backslash \bar{\mathcal{E}}_{l-1}$. Then, $\mathcal{V}^{n_{l}}_C \bigcap \bar{\mathcal{V}}_{l-1} \neq \varnothing$ and $\mathcal{E}^{n_{l}}_C \bigcap \bar{\mathcal{E}}_{l-1} = \varnothing$. By Definition \ref{cycle fixable}, it follows from Lemmas \ref{theorem: bearing and angle for C_n} and \ref{theorem: having same R for positively dependent} that there exists $\mathscr{R}_2\in{\rm O}(2)$ such that $b_{ji}(\bar{q}_{l}) = \mathscr{R}_2b_{ji}(\bar{p}_{l})$ for $(j,i)\in \mathcal{E}^{n_{l}}_C$. It suffices to show that $\mathscr{R}_1=\mathscr{R}_2$. For any $j\in \mathcal{V}^{n_{l}}_C \bigcap \bar{\mathcal{V}}_{l-1}$, there exists $k_1, k_2 \in \mathcal{N}^j_{\bar{\mathcal{G}}_{l-1}}$ such that $b_{jk_1}(\bar{p}_{l})$ is not parallel to $b_{jk_2}(\bar{p}_{l})$. Without loss of generality, let $(k_3,j), (k_4,j)\in \mathcal{E}^{n_{l}}_C$ and $k_1 < k_2 < k_3 < k_4$. Since $(k_1,j,k_3), (k_2,j,k_3), (k_1,j,k_4), (k_2,j,k_4) \in \mathcal{T}_{\bar{\mathcal{G}}_{l+1}}$, one has
    \begin{eqnarray}
       b_{jk_3}(\bar{q}_{l}) &=& \mathscr{R}_1b_{jk_3}(\bar{p}_{l}) = \mathscr{R}_2b_{jk_3}(\bar{p}_{l}), \label{equation: b_{jk_2}} \\
       b_{jk_4}(\bar{q}_{l}) &=& \mathscr{R}_1b_{jk_4}(\bar{p}_{l}) = \mathscr{R}_2b_{jk_4}(\bar{p}_{l}). \label{equation: b_{jk_3}}
    \end{eqnarray}
    According to \cite[Lemma 8.3]{jing2019angle}, (\ref{equation: b_{jk_2}}) implies that $\mathscr{R}_1=\mathscr{R}_2$ or $\mathscr{R}_1\mathscr{H}_{\bar{b}}=\mathscr{R}_2$, where $\bar{b} \triangleq b^\bot_{jk_3}(\bar{p}_{l})$. Suppose that $\mathscr{R}_1\mathscr{H}_{\bar{b}}=\mathscr{R}_2$. Together with (\ref{equation: b_{jk_3}}), we have $b_{jk_4}(\bar{p}_{l}) = \mathscr{H}_{\bar{b}}b_{jk_4}(\bar{p}_{l})$, which holds if and only if $b_{jk_3}(\bar{p}_{l})$ is parallel to $b_{jk_4}(\bar{p}_{l})$ and conflicts with the strict convexity of $(\mathcal{C}^{l}_{n_{l}}, ^l\hspace{-1mm}\bar{p})$.  As a result, $\mathscr{R}_1=\mathscr{R}_2$.

    Case 4b. Suppose that there is a vertex $j \in \mathcal{V}^{n_{l}}_C \bigcap \bar{\mathcal{V}}_{l-1}$ having $(j,i) \in \mathcal{E}^{n_{l}}_C \bigcap \bar{\mathcal{E}}_{l-1}$ and $(j,k) \in \mathcal{E}^{n_{l}}_C \backslash \bar{\mathcal{E}}_{l-1}$. By Definition \ref{cycle fixable}, $(\mathcal{C}^{l}_{n_{l}}, ^l\hspace{-1mm}\bar{p})$ is strictly convex, and there exists $\mathscr{R}_2\in{\rm O}(2)$ such that $b_{ji}(\bar{q}_{l}) = \mathscr{R}_2b_{ji}(\bar{p}_{l})$ for $(j,i)\in \mathcal{E}^{n_{l}}_C$ according to Lemma \ref{theorem: bearing and angle for C_n}. It suffices to show that $\mathscr{R}_1=\mathscr{R}_2$. Since there is a vertex $j$ with $(j,i) \in \mathcal{E}^{n_{l}}_C \bigcap \bar{\mathcal{E}}_{l-1}$ and $(j,k) \in \mathcal{E}^{n_{l}}_C \backslash \bar{\mathcal{E}}_{l-1}$, $b_{ji}(\bar{q}_{l}) = \mathscr{R}_1b_{ji}(\bar{p}_{l}) = \mathscr{R}_2b_{ji}(\bar{p}_{l})$. According to \cite[Lemma 8.3]{jing2019angle}, $\mathscr{R}_1=\mathscr{R}_2$ or $\mathscr{R}_1\mathscr{H}_{\bar{b}}=\mathscr{R}_2$, where $\bar{b} \triangleq b^\bot_{ji}(\bar{p}_{l})$. It follows from Definition \ref{cycle fixable} that there is a vertex $k_1\in \mathcal{N}^j_{\bar{\mathcal{G}}_{l-1}}$ such that $b_{ji}(\bar{p}_{l})$ is not parallel to $b_{jk_1}(\bar{p}_{l})$. Without loss of generality, let $i < k < k_1$, and suppose that $\mathscr{R}_1\mathscr{H}_{\bar{b}}=\mathscr{R}_2$. Then, one has $b_{jk}(\bar{q}_{l}) = \mathscr{R}_2b_{jk}(\bar{p}_{l}) = \mathscr{R}_1 \mathscr{H}_{\bar{b}}b_{jk}(\bar{p}_{l})$. From $b_{jk_1}(\bar{q}_{l}) = \mathscr{R}_1b_{jk_1}(\bar{p}_{l})$ and $(k,j,k_1)\in\mathcal{T}_{\mathcal{G}}$,
    \begin{eqnarray}
       b^\top_{jk_1}(\bar{q}_{l})b_{jk}(\bar{q}_{l}) = b^\top_{jk_1}(\bar{p}_{l})\mathscr{H}_{\bar{b}}b_{jk}(\bar{p}_{l}) = b^\top_{jk_1}(\bar{p}_{l})b_{jk}(\bar{p}_{l}),
    \end{eqnarray}
    which implies $b_{ji}(\bar{p}_{l})$ is parallel to $b_{jk_1}(\bar{p}_{l})$ or $b_{jk}(\bar{p}_{l})$, contradicting Definition \ref{cycle fixable}.
    
    Case 5. $|\mathcal{V}^{n_{l}}_C \bigcap \bar{\mathcal{V}}_{l-1}| = 1$. The proof is similar to that of Case 4a.
\end{proof}

\section{Proof of Theorem \ref{theorem: equivalence in framework}} \label{proof of theorm: theorem: equivalence in framework} 
\begin{proof}
    The necessity has been proven in \cite[Theorem 5]{jing2019angle}. It suffices to show that $(\mathcal{G},p)$ is (GAR) AR if it is (GBR) BR.
         
        Suppose $(\mathcal{G},p)$ is BR. Then, there exists a neighborhood $U_p$ of $p$ such that $B^{-1}_{\mathcal{G}}(B_{\mathcal{G}}(p)) \bigcap U_p = B^{-1}_{\mathcal{K}}(B_{\mathcal{K}}(p)) \bigcap U_p$. Note that $A^{-1}_{\mathcal{K}}(A_{\mathcal{K}}(p)) \subseteq A^{-1}_{\mathcal{G}}(A_{\mathcal{G}}(p))$. Therefore, it suffices to show that for this $U_p$, $A^{-1}_{\mathcal{G}}(A_{\mathcal{G}}(p))\bigcap U_p \subseteq A^{-1}_{\mathcal{K}}(A_{\mathcal{K}}(p))\bigcap U_p$. For any $q \in A^{-1}_{\mathcal{G}}(A_{\mathcal{G}}(p))\bigcap U_p$, since $(\mathcal{G}, p)$ is cycle compatible, it follows from Lemma \ref{lemma: cycle fixable can determine the bearing} that $(I_{m}\otimes\mathscr{R})^{-1}q \in B^{-1}_{\mathcal{G}}(B_{\mathcal{G}}(p)) \bigcap U_p$ for some $\mathscr{R} \in {\rm O}(2)$. Since $B^{-1}_{\mathcal{G}}(B_{\mathcal{G}}(p)) \bigcap U_p= B^{-1}_{\mathcal{K}}(B_{\mathcal{K}}(p))\bigcap U_p$, one has $(I_n\otimes\mathscr{R})^{-1}q \in B^{-1}_{\mathcal{K}}(B_{\mathcal{K}}(p)) \bigcap U_p$, i.e., $b_{ji}(q) = \mathscr{R}b_{ji}(p)$ for any $j,i\in\mathcal{V}$. Thus, for any $(i,j,k)\in\mathcal{T}_{\mathcal{K}}$, $b^\top_{ji}(q)b_{jk}(q) = b^\top_{ji}(p)\mathscr{R}^\top\mathscr{R}b_{jk}(p) = b^\top_{ji}(p)b_{jk}(p)$, which implies that $q\in A^{-1}_{\mathcal{K}}(A_{\mathcal{K}}(p))\bigcap U_p$. 
    
    By extending $U_p$ to $\mathbb{R}^{2n}$, the same result can be verified for GAR and GBR.
\end{proof}

\section{Proof of Lemma \ref{lemma: angle can determine the bearing}} \label{poof of lemma: angle can determine the bearing} 
    \begin{proof}
    Sufficiency follows from $B^{-1}_{\mathcal{G}}(B_{\mathcal{G}}(p)) \subseteq A^{-1}_{\mathcal{G}}(A_{\mathcal{G}}(p))$.

    Necessity. We first show that for each $j\in\mathcal{V}$, one has $b_{ji}(q) = \mathscr{R}_jb_{ji}(p)$ for $i\in\mathcal{N}_j$ and some $\mathscr{R}_j\in{\rm O}(2)$. Since $p\in\mathbb{R}^{2n}$ is generic, there exist two neighbors $i,k\in\mathcal{N}_j$ of $j$ with $b^\top_{ji}(q)b_{jk}(q) = b^\top_{ji}(p)b_{jk}(p) \neq \pm 1$. Let $G_{j} = (\dots, b_{jh}, \dots)\in \mathbb{R}^{2\times |\mathcal{N}_j|}$ for $h\in \mathcal{N}_j$. Since $A_{\mathcal{G}}(q) = A_{\mathcal{G}}(p)$ and $(i,j,k)\in\mathcal{T}_{\mathcal{G}}$ for any pair of edges $(j,i)$ and $(j,k)$, we have $G_{j}^\top(q)G_{j}(q) = G_{j}^\top(p)G_{j}(p)$. It follows from $\rank(G_{j}(p)) = 2$ and \cite[Lemma 8.1]{jing2019angle} that $G_{j}(q) = \mathscr{R}_jG_{j}(p)$ for some $\mathscr{R}_j\in{\rm O}(2)$. That is, $b_{ji}(q) = \mathscr{R}_jb_{ji}(p)$ for $i\in\mathcal{N}_j$. Hence, it suffices to show that for any $(i,j)\in\mathcal{E}$, one has $\mathscr{R}_i = \mathscr{R}_j$. Next, we show that if $\mathscr{R}_i \neq \mathscr{R}_j$ for some $(i,j) \in \mathcal{E}$, there is a non-zero polynomial $h(x_1,\dots,x_n)$ with integer coefficients such that $h(p_1,\dots,p_n) = 0$, which conflicts with the fact that $p$ is a generic configuration.

    We first note that the vector $\mathscr{R}_jz_{ji} \triangleq\mathscr{R}_j(p_i-p_j) \in\mathbb{R}^2$ can be described in the complex plane $\mathbb{C}$ uniquely, i.e.,
	\begin{eqnarray}
		\label{eq:complex transformation}
		\mathscr{R}_jz_{ji} & =\left\{\begin{array}{ll}
			a_jc_{ji} & {\rm if}\hspace{2mm}{\rm det}(\mathscr{R}_j) = 1, \\
			a_j\bar{c}_{ji} & {\rm if}\hspace{2mm}{\rm det}(\mathscr{R}_j) = -1,
		\end{array}\right.
	\end{eqnarray} 
    where $a_j = e^{\theta_j{\rm i}}$ for ${\rm i} = \sqrt{-1}$ and some $\theta_j\in(0,2\pi]$, $\bar{c}_{ji}$ is the conjugate of $c_{ji}$, and $c_{ji} = (x_i-x_j)+(y_i-y_j){\rm i}$. Considering an edge $(i,j)\in\mathcal{E}$, one has $b_{ji}(q) = \mathscr{R}_jb_{ji}(p) = \mathscr{R}_ib_{ji}(p)$, and hence, $\mathscr{R}_jz_{ji} = \mathscr{R}_iz_{ji}$. Describing $\mathscr{R}_jz_{ji} = \mathscr{R}_iz_{ji}$ in the complex plane by (\ref{eq:complex transformation}), there may exist three cases:
        \begin{itemize}
            \item Case 1. ${\rm det}(\mathscr{R}_j) = {\rm det}(\mathscr{R}_i)$, one has $\frac{c_{ji}}{c_{ji}} = \frac{a_{i}}{a_{j}} = 1$ or $\frac{\bar{c}_{ji}}{\bar{c}_{ji}} = \frac{a_{i}}{a_{j}} = 1$;
            \item Case 2. ${\rm det}(\mathscr{R}_j) = 1$ and ${\rm det}(\mathscr{R}_i) = -1$, one has $\frac{c_{ji}}{\bar{c}_{ji}} = \frac{a_{i}}{a_{j}}$;
            \item Case 3. ${\rm det}(\mathscr{R}_j) = -1$ and ${\rm det}(\mathscr{R}_i) = 1$, one has $\frac{\bar{c}_{ji}}{c_{ji}} = \frac{a_{i}}{a_{j}}$.
        \end{itemize}

    Suppose that there exist edges $(i,j)\in\mathcal{E}$ satisfying $\mathscr{R}_i \neq \mathscr{R}_j$. Since $\mathcal{G}$ is a 2-edge-connected graph, these edges must belong to some cycle frameworks of $(\mathcal{G},q)$. Without loss of generality, let $(\mathcal{C}_s,\bar{q})$ be a cycle subframework of the framework $(\mathcal{G},q)$ such that there exists $(i,j)\in\mathcal{E}^s_C$ satisfying $\mathscr{R}_i \neq \mathscr{R}_j$. For each edge $(j,i)\in\mathcal{E}^s_C$, one has $\mathscr{R}_jz_{ji} = \mathscr{R}_iz_{ji}$, and this equation can be written as one of the three equations mentioned above in the complex plane. Multiplying all the equations for $(j,i)\in\mathcal{E}^s_C$ together gives the product equation
    \begin{eqnarray} \label{equation: cycle and R}
        \prod \limits_{(i,j)\in\mathcal{E}_+}\left(\frac{c_{ji}}{\bar{c}_{ji}}\right)\prod \limits_{(i,j)\in\mathcal{E}_-}\left(\frac{\bar{c}_{ji}}{c_{ji}}\right) &=& \frac{a_{2}}{a_{1}}\frac{a_{3}}{a_{2}}\cdots \frac{a_{s}}{a_{s-1}}\frac{a_{1}}{a_{s}} = 1.
    \end{eqnarray}
    where $\mathcal{E}_+ = \left\{(j,i)\in \mathcal{E}^s_C: {\rm det}(\mathscr{R}_j) = 1, {\rm det}(\mathscr{R}_i) = -1\right\}$ and $\mathcal{E}_- = \{(j,i)\in \mathcal{E}^s_C: {\rm det} (\mathscr{R}_j) \\  = -1, {\rm det}(\mathscr{R}_i) = 1\}$. Note that both $c_{ji}$ and $\bar{c}_{ji}$ are irreducible polynomials related to variables $p_i$ and $p_j$. Similarly, since each edge appears only once in (\ref{equation: cycle and R}), $c_{ji}$ (resp. $\bar{c}_{ji}$) and $\bar{c}_{gh}$ (resp. $c_{gh}$) are irreducible for distinct edges $(j,i), (g,h)\in\mathcal{E}^s_C$ as well. As a result, (\ref{equation: cycle and R}) is not an identity equation. By performing algebraic manipulation on \eqref{equation: cycle and R}, we eliminate the denominators and obtain
    \begin{eqnarray} \label{equation: cycle and R2}
        \prod \limits_{(i,j)\in\mathcal{E}_+}c_{ji}\prod \limits_{(i,j)\in\mathcal{E}_-}\bar{c}_{ji}-\prod \limits_{(i,j)\in\mathcal{E}_+}\bar{c}_{ji}\prod \limits_{(i,j)\in\mathcal{E}_-}c_{ji} = 0.
    \end{eqnarray}
    Let $H(p)$ be the left side of (\ref{equation: cycle and R2}), then there exists a non-zero polynomial $h(x_1,\dots,x_n)$ with integer coefficients satisfying $H(p) = 2h(p){\rm i} = 0$, which implies that $h(p) = 0$ and conflicts with the fact that $p$ is a generic configuration.
    \end{proof}

\section{Notation of graph theory} \label{Notation of graph theory}Following from \cite{murty2008graph}, a \textit{plane graph} $\tilde{\mathcal{G}} = (\tilde{\mathcal{V}}, \tilde{\mathcal{E}})$ is a planar graph $\mathcal{G}$ together with a fixed planar embedding in $\mathbb{R}^2$. Fig. \ref{illustration_of_algorithm_1} shows plane Laman graphs. A plane graph $\tilde{\mathcal{G}}$ partitions the plane into a number of arc-wise-connected open sets by the simple curves and points. These sets are called the \textit{faces} $\mathcal{F}$ of $\tilde{\mathcal{G}}$, and can be classified into the unbounded outer face $\mathcal{F}_{\rm out}$ and bounded inner faces $\mathcal{F}_{\rm in}$. In Fig. \ref{illustration_of_algorithm_1}, $\mathcal{F}_1$ is the outer face, and the rest are inner faces.

\begin{definition}
    Given a plane graph $\tilde{\mathcal{G}}$, the boundary $\partial(\mathcal{F}_{\rm in})$ of an inner face $\mathcal{F}_{\rm in}$ is the boundary of the associated bounded open set, while the boundary $\partial(\mathcal{F}_{\rm out})$ of the outer face $\mathcal{F}_{\rm out}$ is the boundary of the closed set $\mathcal{O}$ whose complement set relative to $\mathbb{R}^2$ is the unbounded open set associated with $\mathcal{F}_{\rm out}$.
\end{definition}

When there is no scope for confusion, we use the boundary $\partial(\mathcal{F})$ to denote the corresponding cycle $\mathcal{C}^l_{n_l}$ of the two-edge-connected graph $\mathcal{G}$. Accordingly, the length of $\partial(\mathcal{F})$ is denoted by $|\partial(\mathcal{F})| = n_l$. For example, the boundary of face $\mathcal{F}_1$ in Fig. \ref{illustration_of_algorithm_1}(a) corresponds to the cycle composed of vertices $1,2$, and 3, of which the length is 3. On the other hand, two faces are adjacent if their boundaries share an edge $(i,j)\in\mathcal{E}$. 

\begin{definition}[\cite{hedetniemi1985interior}]
    The dual graph $\hat{\mathcal{G}}(\tilde{\mathcal{G}}) = (\mathcal{V}(\tilde{\mathcal{G}}),\mathcal{E}(\tilde{\mathcal{G}}))$ of a plane graph $\tilde{\mathcal{G}}$ is an undirected graph, where the elements $\mathcal{F}_i \in \mathcal{V}(\tilde{\mathcal{G}})$ correspond to the faces of $\tilde{\mathcal{G}}$ and the elements $(\mathcal{F}_i,\mathcal{F}_j)\in \mathcal{E}(\tilde{\mathcal{G}})$ correspond to the adjacency of the faces in $\tilde{\mathcal{G}}$.
\end{definition}

Fig. \ref{illustration_of_algorithm_1} provides some examples of the dual graph. The \textit{weak dual graph} $\hat{\mathcal{G}}(\tilde{\mathcal{G}})$ of a plane graph is obtained from its dual graph by deleting the vertex corresponding to the outer face, together with all edges incident to this vertex. If $\tilde{\mathcal{G}}$ is clear from the context, we simply denote $\hat{\mathcal{G}}(\tilde{\mathcal{G}}) = (\mathcal{V}(\tilde{\mathcal{G}}),\mathcal{E}(\tilde{\mathcal{G}}))$ by $\hat{\mathcal{G}} = (\hat{\mathcal{V}},\hat{\mathcal{E}})$. 

\section{Proof of Theorem \ref{theorem: proof for algorithm}} \label{proof of theorem: proof for algorithm} 
\begin{proof}
    It suffices to show that the weak dual graph $\hat{\mathcal{G}}$ of the plane graph $\tilde{\mathcal{G}}$ produced by Algorithm \ref{alg:1} is connected, and $\hat{\mathcal{G}}$ always satisfies $|\partial(\mathcal{F}_i)| \leq 4$ for $\mathcal{F}_i\in \hat{\mathcal{V}}$. Then, we can induce a set of cycles satisfying $\mathcal{L}_h = \bigcup^c_{l=1}\mathcal{C}^{l}_{n_l}$ by finding a tree in the weak dual graph $\hat{\mathcal{G}}$. We prove this by induction.

    For the initial plane graph with $h=3$ vertices, the statement holds directly. 
    
    Suppose that the weak dual graph $\hat{\mathcal{G}}$ of the plane graph $\tilde{\mathcal{G}}$ produced by Algorithm \ref{alg:1} with $h\geq3$ is connected, and $|\partial(\mathcal{F}_i)| \leq 4$ for $\mathcal{F}_i\in \hat{\mathcal{V}}$. Next, we show that by introducing a new vertex $v$ to $\tilde{\mathcal{G}}$ via either vertex addition or edge splitting, the updated plane graph $\tilde{\mathcal{G}}^\prime$ preserves the same property.
    
    Case 1. Performing the vertex addition operation splits the outer face $\mathcal{F}_{\rm out}$ of $\tilde{\mathcal{G}}$ into two faces: a new outer face and an inner face $\mathcal{F}^\prime$. $\partial(\mathcal{F}^\prime)$ shares at least one edge with $\partial(\mathcal{F}_i)$ for some inner face $\mathcal{F}_i\in\hat{\mathcal{V}}$. Together with $|\partial(\mathcal{F}^\prime)| \leq 4$, the result is preserved for the updated plane graph $\tilde{\mathcal{G}}^\prime$.

    Case 2. Performing the edge splitting operation divides the merged face of $\mathcal{F}_1$ and $\mathcal{F}_2$ into three new faces $\mathcal{F}^\prime_i, i\in\{1,2,3\}$, any two of which share an edge. Since $\mathcal{F}^\prime_1$ with $|\partial(\mathcal{F}^\prime_1)| = 4$ is generated from the selected face $\mathcal{F}_1$ with $|\partial(\mathcal{F}_1)| = 3$ in $\tilde{\mathcal{G}}$, the weak dual graph of the updated plane graph $\tilde{\mathcal{G}}^\prime$ is connected. Since $|\partial(\mathcal{F}_i)| \leq 4$ for $\mathcal{F}_i\in \hat{\mathcal{V}}$, there are two cases for $\mathcal{F}_1$ and $\mathcal{F}_2$: (\rmnum{1}) $\partial(\mathcal{F}_1)$ shares one edge with $\partial(\mathcal{F}_2)$ and $|\partial(\mathcal{F}_2)| = 3$; (\rmnum{2}) $\partial(\mathcal{F}_1)$ shares at most two edges with $\partial(\mathcal{F}_2)$ where $\mathcal{F}_2  = \mathcal{F}_{\rm out}$ or $|\partial(\mathcal{F}_2)| = 4$. Note that $\mathcal{F}^\prime_2$ and $\mathcal{F}^\prime_3$ consist of the rest part of the merged face, where $\partial(\mathcal{F}^\prime_2)$ satisfies $|\partial(\mathcal{F}^\prime_2)| \leq 4$ due to the existence of a path from $i$ to $k$ of length at most $2$ along $\partial(\mathcal{F}_2)$. It can be verified that in these cases, after performing the edge splitting operation, $\mathcal{F}^\prime_3$ satisfies $\mathcal{F}^\prime_3 = \mathcal{F}_{\rm out}$ or $|\partial(\mathcal{F}^\prime_3)| \leq 4$ in $\tilde{\mathcal{G}}^\prime$.   
\end{proof}

\section{Proof of Theorem \ref{theorem:localization theorem}} \label{proof of theorem:localization theorem}
\begin{proof}

From \cite[Corollary 3]{zhao2016localizability}, Assumption \ref{ass: 1} implies that $(\mathcal{G},p,\mathcal{A})$ is bearing localizable. Since the estimator (\ref{eq:location estimators2}) is independent of the position estimates $\hat{p}_j(t)$ for all $j\in\mathcal{F}$ and $t\geq0$, it suffices to show that $\hat{b}_{ji}(t)$ for any $(i,j) \in \mathcal{E}\backslash\mathcal{E}^{n_1}_C$ asymptotically converges to the global bearing $b_{ji}$ with $t \rightarrow\infty$, and the convergence of $\hat{p}_j$ to $p_j$ for all $j\in \mathcal{F}$ then follows by an argument similar to that in \cite[Theorem 2]{li2019globally}.

We first note that since the root edges $(r_1, k_1)$ and $(r_2, k_2)$ of $(i,j) \in \mathcal{E}\backslash\mathcal{E}^{n_1}_C$ are not parallel, $s_{ji} \triangleq {\rm sin}\beta_{ji} = {\rm sin}(\alpha_{r_2k_2}-\alpha_{r_1k_1}) \neq 0$ and one has 
    \begin{eqnarray*}
        b_{ji} &=& \begin{bmatrix}
        b^\top_{r_1k_1} \\
        b^\top_{r_2k_2} \\
        \end{bmatrix}^{-1}\begin{bmatrix}
        c_1 \\
        c_2 \\
        \end{bmatrix} =  \begin{bmatrix}
        {\rm cos}\alpha_{r_1k_1} & {\rm sin}\alpha_{r_1k_1} \\
        {\rm cos}\alpha_{r_2k_2} & {\rm sin}\alpha_{r_2k_2} \\
        \end{bmatrix}^{-1}\begin{bmatrix}
        c_1 \\
        c_2 \\
        \end{bmatrix} = \frac{1}{s_{ji}} \begin{bmatrix}
        -b^{\bot}_{r_2k_2} & b^{\bot}_{r_1k_1} \\
        \end{bmatrix}\begin{bmatrix}
        c_1 \\
        c_2 \\
        \end{bmatrix}.  
    \end{eqnarray*}
where $c_1 \triangleq {\rm cos}\theta_{r_1k_1,ji}$ and $c_2 \triangleq {\rm cos}\theta_{r_2k_2,ji}$. Multiplying both sides by $s_{ji}$, one has
    \begin{eqnarray} \label{equation: estimation equation}
        s_{ji}b_{ji} = \mathscr{R}_o(\frac{\pi}{2})(c_2b_{r_1k_1}-c_1b_{r_2k_2}),
    \end{eqnarray}
Let $\delta_{ji} \triangleq \hat{b}_{ji} - b_{ji}$ and $\hat{s}_{ji} \triangleq {\rm sin}\hat{\beta}_{ji}$. It follows from (\ref{eq:location estimators2}) and (\ref{equation: estimation equation}) that 
    \begin{eqnarray}
        \dot{\delta}_{ji} = -\hat{s}_{ji}^2\delta_{ji} + w(t). \label{eq:proof of the estimator 7}
    \end{eqnarray}
where $w(t) \triangleq \hat{s}_{ji}(s_{ji}-\hat{s}_{ji})b_{ji} + \hat{s}_{ji}\mathscr{R}_o(\frac{\pi}{2})(c_2\delta_{r_1k_1} - c_1\delta_{r_2k_2})$. Next, we show that $\delta_{ji}$ asymptotically converges to zero by induction.

    For $\{r_1,k_1,r_2,k_2\} \subseteq \mathcal{A}$, it follows that $s_{ji}=\hat{s}_{ji} \neq 0$ and $\delta_{r_1k_1} = \delta_{r_2k_2} = 0$. Thus, according to (\ref{eq:proof of the estimator 7}), $\dot{\delta}_{ji} = -s_{ji}^2\delta_{ji}$ and $\delta_{ji}$ converges to zero exponentially fast, which also implies that $\hat{b}_{ji}$ is continuous in $[0,\infty)$.

    Suppose that for $\{r_1,k_1,r_2,k_2\} \backslash \mathcal{A} \neq \varnothing$, $\hat{b}_{r_1k_1}$ and $\hat{b}_{r_2k_2}$ are continuous in $[0, H]$ for any $0< H < \infty$ and the estimated errors of the root edges asymptotically converge to zero. Next, we show that $\delta_{ji}$ is asymptotically converge to zero.
    
    To this end, we first show that $\hat{b}_{ji}$ is continuous in $[0, H]$. Since $\hat{b}_{r_1k_1}$ and $\hat{b}_{r_2k_2}$ are continuous in $[0, H]$ and $\|\hat{s}_{ji}\|\in[0,1]$, $r(t) \triangleq \hat{s}_{ji}\mathscr{R}_o(\frac{\pi}{2})(c_2\hat{b}_{r_1k_1} - c_1\hat{b}_{r_2k_2})$ and $\hat{s}_{ji}$ are both bounded and measurable in $[0, H]$. It follows from (\ref{eq:location estimators2}) that, 
    \begin{eqnarray}
\hat b_{ji}(t) = e^{-A(t)} \left[\hat b(0)+\underbrace{\int_0^t e^{A(\tau)}r(\tau)\,\mathrm{d}\tau}_{\chi_{ji}(t)}\right],
    \end{eqnarray}
    where $t\in[0,H]$ and $A(t) \triangleq \int_0^t \hat{s}^2_{ji}(\xi)\,\mathrm{d}\xi$. Note that $A(t)$ is Lipschitz continuous and satisfies $A(t) \in [0,H]$. Thus, $e^{A(t)} r(t)$ is bounded. Hence, $\chi_{ji}(t)$ is Lipschitz continuous, and it can be verified that $\hat b_{ji}$ is also Lipschitz continuous in $[0, H]$.

    Next, we show the convergence. There exists a finite time $T>0$ such that for all $t>T$, it always holds that $\hat{s}^2_{ji}(t) \geq s_0 \triangleq s^2_{ji}/2$. Since $w(t)$ is bounded and measurable in $[0, H]$ for any $0< H < \infty$, together with (\ref{eq:proof of the estimator 7}), it holds that for any $t\geq T$, one has 
    \begin{eqnarray}
    \delta_{ji}(t)= e^{-\int_T^t a(\xi)\,d\xi}\delta_{ji}(T) + \int_T^te^{-\int_\tau^t a(\xi)\,d\xi}w(\tau)\,d\tau,
    \end{eqnarray}
    where $\delta_{ji}(T)$ is bounded since $\hat{b}_{ji}$ is continuous in $[0, T]$. Note that $e^{-\int_\tau^t a(\xi)\,d\xi} \leq e^{-s_0(t-\tau)}$. Hence, one has
    \begin{eqnarray*}
    \left\|\int_T^te^{-\int_\tau^t a(\xi)\,d\xi}w(\tau)\,d\tau \right\| \leq \int_T^te^{-\int_\tau^t a(\xi)\,d\xi}\|w(\tau)\|d\tau \leq \frac{1}{s_0} \sup_{T\leq \tau\leq t}\|w(\tau)\|, 
    \end{eqnarray*}
    As a result, 
    \begin{eqnarray}
    \|\delta_{ji}(t)\| \leq e^{-s_0(t-T)}\|\delta_{ji}(T)\| + \frac{1}{s_0}\sup_{T\leq \tau\leq t}\|w(\tau)\|.
    \end{eqnarray}
     Since $w(t) \rightarrow 0$ as $t \rightarrow \infty$, according to input-to-state stable \cite[pp. 175]{khalil2002nonlinear}, it can be concluded that $\delta_{ji}$ converges to zero asymptotically.
\end{proof}
\section*{Acknowledgments}
The authors would like to thank Yongjie Liu for the valuable discussions on the connection between bearing constraints and angle constraints within convex cycle frameworks.

\bibliographystyle{siamplain}
\bibliography{references}
\end{document}

%% file: ex_shared.tex
\usepackage{lipsum}
\usepackage{epstopdf}

\usepackage{mathrsfs} 
\usepackage{amsmath}
\usepackage{amsfonts, amssymb}
\usepackage{graphicx}
\usepackage{color}
\usepackage{algorithm}
\usepackage{algorithmic}

\newtheorem{problem}{Problem}
\newtheorem{assumption}{Assumption}

\def\BibTeX{{\rm B\kern-.05em{\sc i\kern-.025em b}\kern-.08em
		T\kern-.1667em\lower.7ex\hbox{E}\kern-.125emX}}
\allowdisplaybreaks
\DeclareMathOperator{\rank}{rank}

\usepackage[numbers, sort&compress]{natbib}

\usepackage{multirow}
\usepackage{nicematrix}

\makeatletter
\newcommand{\rmnum}[1]{\romannumeral #1}
\newcommand{\Rmnum}[1]{\expandafter\@slowromancap\romannumeral #1@}
\makeatother

\ifpdf
  \DeclareGraphicsExtensions{.eps,.pdf,.png,.jpg}
\else
  \DeclareGraphicsExtensions{.eps}
\fi

\newsiamremark{remark}{Remark}
\newsiamremark{hypothesis}{Hypothesis}
\crefname{hypothesis}{Hypothesis}{Hypotheses}
\newsiamthm{claim}{Claim}
\newsiamremark{fact}{Fact}
\crefname{fact}{Fact}{Facts}

\headers{Equivalent Graphical Conditions for Angle Rigidity}{J. HUANG, G. JING, K. ZHAO, and Y. SONG}

\title{Equivalent Graphical Conditions for Angle Rigidity and Angle-Based Network Localization\thanks{Submitted to the editors DATE.
\funding{This work was supported in part by the National Key Research and Development Program of China under Grant 2025YFA1018800 and in part by the National Natural Science Foundation of China under Grant 62573068 and Grant 62533002.}}}

\author{JINPENG HUANG\thanks{School of Automation, Chongqing University, Chongqing 400044, China
  (\email{202313021067t@stu.cqu.edu.cn}, \email{zhaokai@cqu.edu.cn}, \email{ydsong@cqu.edu.cn}).}
\and GANGSHAN JING\thanks{Corresponding author. School of Automation, Chongqing University, Chongqing 400044, China 
  (\email{jinggangshan@cqu.edu.cn}).}
\and KAI ZHAO\footnotemark[2]
  \and YONGDUAN SONG\footnotemark[2]}

\usepackage{amsopn}


%% file: main.bbl
\begin{thebibliography}{10}

\bibitem{anderson2008rigid}
{\sc B.~D.~O. Anderson, C.~Yu, B.~Fidan, and J.~M. Hendrickx}, {\em Rigid graph
  control architectures for autonomous formations}, IEEE Control Syst. Mag., 28
  (2008), pp.~48--63.

\bibitem{arrigoni2018bearing}
{\sc F.~Arrigoni and F.~Andrea}, {\em Bearing-based network localizability: A
  unifying view}, IEEE Trans. Pattern Anal. Mach. Intell., 41 (2018),
  pp.~2049--2069.

\bibitem{asimow1978rigidity}
{\sc L.~Asimow and B.~Roth}, {\em The rigidity of graphs}, Trans. Amer. Math.
  Soc., 245 (1978), pp.~279--289.

\bibitem{aspnes2006theory}
{\sc J.~Aspnes, T.~Eren, D.~K. Goldenberg, A.~S. Morse, W.~Whiteley, Y.~R.
  Yang, B.~D.~O. Anderson, and P.~N. Belhumeur}, {\em A theory of network
  localization}, IEEE Trans. Mob. Comput., 5 (2006), pp.~1663--1678.

\bibitem{bickle2020fundamentals}
{\sc A.~Bickle}, {\em Fundamentals of graph theory}, vol.~43, American
  Mathematical Soc., 2020.

\bibitem{murty2008graph}
{\sc J.~Bondy and A.~Murty}, {\em Graph Theory {\rm (Graduate Texts in
  Mathematics 244)}}, vol.~244, Springer, 2008.

\bibitem{cao2019ratio}
{\sc K.~Cao, Z.~Han, X.~Li, and L.~Xie}, {\em Ratio-of-distance rigidity theory
  with application to similar formation control}, IEEE Trans. Autom. Control,
  65 (2019), pp.~2598--2611.

\bibitem{chen2022triangular}
{\sc L.~Chen}, {\em Triangular angle rigidity for distributed localization in
  {2D}}, Automatica, 143 (2022), p.~110414.

\bibitem{chen2020angle}
{\sc L.~Chen, M.~Cao, and C.~Li}, {\em Angle rigidity and its usage to
  stabilize multiagent formations in {2-D}}, IEEE Trans. Autom. Control, 66
  (2020), pp.~3667--3681.

\bibitem{eren2012formation}
{\sc T.~Eren}, {\em Formation shape control based on bearing rigidity}, Int. J.
  Control, 85 (2012), pp.~1361--1379.

\bibitem{fang2009sequential}
{\sc J.~Fang, M.~Cao, A.~S. Morse, and B.~D.~O. Anderson}, {\em Sequential
  localization of sensor networks}, SIAM J. Control Optim., 48 (2009),
  pp.~321--350.

\bibitem{fang20232}
{\sc X.~Fang, J.~Li, X.~Li, and L.~Xie}, {\em 2-d distributed pose estimation
  of multi-agent systems using bearing measurements}, J. Autom. Intell., 2
  (2023), pp.~70--78.

\bibitem{hedetniemi1985interior}
{\sc S.~M. Hedetniemi, A.~Proskurowski, and M.~M. Sys{\l}o}, {\em Interior
  graphs of maximal outerplane graphs}, Journal of Combinatorial Theory, Series
  B, 38 (1985), pp.~156--167.

\bibitem{huang2026signed}
{\sc J.~Huang and G.~Jing}, {\em Signed angle rigid graphs for network
  localization and formation control}, IEEE Trans. Autom. Control, 71 (2026),
  pp.~6369--6384.

\bibitem{jackson2005connected}
{\sc B.~Jackson and T.~Jord{\'a}n}, {\em Connected rigidity matroids and unique
  realizations of graphs}, J. Combin. Theory Ser. B, 94 (2005), pp.~1--29.

\bibitem{jing2021angle}
{\sc G.~Jing, C.~Wan, and R.~Dai}, {\em Angle-based sensor network
  localization}, IEEE Trans. Autom. Control, 67 (2021), pp.~840--855.

\bibitem{jing2018weak}
{\sc G.~Jing, G.~Zhang, H.~W. Lee, and L.~Wang}, {\em Weak rigidity theory and
  its application to formation stabilization}, SIAM J. Control Optim., 56
  (2018), pp.~2248--2273.

\bibitem{jing2019angle}
{\sc G.~Jing, G.~Zhang, H.~W.~J. Lee, and L.~Wang}, {\em Angle-based shape
  determination theory of planar graphs with application to formation
  stabilization}, Automatica, 105 (2019), pp.~117--129.

\bibitem{khalil2002nonlinear}
{\sc H.~K. Khalil}, {\em \emph{Nonlinear systems}, 3rd ed}.
\newblock Englewood Cliffs, NJ, USA: Prentice-Hall, 2002.

\bibitem{ko2020bearing}
{\sc G.~Ko, M.~H. Trinh, and H.~S. Ahn}, {\em Bearing-only control of directed
  cycle formations: almost global convergence and hardware implementation},
  Int. J. Robust Nonlinear Control, 30 (2020), pp.~4789--4804.

\bibitem{krick2009stabilisation}
{\sc L.~Krick, M.~E. Broucke, and B.~A. Francis}, {\em Stabilisation of
  infinitesimally rigid formations of multi-robot networks}, Int. J. Control,
  82 (2009), pp.~423--439.

\bibitem{li2019globally}
{\sc X.~Li, X.~Luo, and S.~Zhao}, {\em Globally convergent distributed network
  localization using locally measured bearings}, IEEE Trans. Control Netw.
  Syst., 7 (2019), pp.~245--253.

\bibitem{liu2025rigidity}
{\sc Y.~Liu, G.~Jing, and L.~Wang}, {\em Rigidity theory of graphs with
  heterogeneous vertices}, arXiv preprint arXiv:2510.27515,  (2025).

\bibitem{michieletto2021unified}
{\sc G.~Michieletto, A.~Cenedese, and D.~Zelazo}, {\em A unified dissertation
  on bearing rigidity theory}, IEEE Trans. Control Netw. Syst., 8 (2021),
  pp.~1624--1636.

\bibitem{oh2026subtended}
{\sc K.-H. Oh, B.~Fidan, J.-G. Lee, and H.-S. Ahn}, {\em Subtended angle-based
  network localization via distributed bearing vector estimation}, Systems
  Control Lett., 215 (2026), p.~106489.

\bibitem{oh2015survey}
{\sc K.~K. Oh, M.~C. Park, and H.~S. Ahn}, {\em A survey of multi-agent
  formation control}, Automatica, 53 (2015), pp.~424--440.

\bibitem{peng2025angle}
{\sc C.~Peng and J.~Huang}, {\em Angle-constrained distributed formation
  control with disturbance rejection}, Automatica, 177 (2025), p.~112314.

\bibitem{trinh2018formations}
{\sc M.~H. Trinh, D.~Mukherjee, D.~Zelazo, and H.-S. Ahn}, {\em Formations on
  directed cycles with bearing-only measurements}, Int. J. Robust Nonlinear
  Control, 28 (2018), pp.~1074--1096.

\bibitem{wan2019sensor}
{\sc C.~Wan, G.~Jing, S.~You, and R.~Dai}, {\em Sensor network localization via
  alternating rank minimization algorithms}, IEEE Trans. Control Netw. Syst., 7
  (2019), pp.~1040--1051.

\bibitem{zhao2017laman}
{\sc S.~Zhao, Z.~Sun, D.~Zelazo, M.~H. Trinh, and H.~S. Ahn}, {\em Laman graphs
  are generically bearing rigid in arbitrary dimensions}, in 2017 IEEE 56th
  Annual Conference on Decision and Control (CDC), 2017, pp.~3356--3361.

\bibitem{zhao2015bearing}
{\sc S.~Zhao and D.~Zelazo}, {\em Bearing rigidity and almost global
  bearing-only formation stabilization}, IEEE Trans. Autom. Control, 61 (2015),
  pp.~1255--1268.

\bibitem{zhao2016localizability}
{\sc S.~Zhao and D.~Zelazo}, {\em Localizability and distributed protocols for
  bearing-based network localization in arbitrary dimensions}, Automatica, 69
  (2016), pp.~334--341.

\end{thebibliography}
